\documentstyle[12pt]{article}   %%% PREPRINT %%%
\newcommand{\be}{\begin{equation}}
\newcommand{\ee}{\end{equation}}
\newcommand{\bea}{\begin{eqnarray}}
\newcommand{\eea}{\end{eqnarray}}
\newcommand{\bean}{\begin{eqnarray*}}
\newcommand{\eean}{\end{eqnarray*}}
\newcommand{\N}{I\!\!N}
\newcommand{\PP}{I\!\!P}
\newcommand{\R}{I\!\!R}
\newcommand{\Z}{Z\!\!\!Z}

\newcommand{\C}{\,I\!\!\!\!C}
\newcommand{\Q}{I\!\!\!\!Q}

\newcommand{\CM}{{\cal M}}

\newcommand{\veee}{\scriptscriptstyle\vee}

\newcounter{Abschnitt}[section]

\newcommand{\neu}[1]{\protect\refstepcounter{Abschnitt}\protect\label{t#1}\vspace{1ex}
                     {\bf (\arabic{section}.\protect\arabic{Abschnitt})}
                     $\qquad$}
\newcommand{\zitat}[2]{(\protect\ref{s#1}.\protect\ref{t#1#2})}
\newcommand{\surj}{\longrightarrow\hspace{-1.5em}\longrightarrow}
\newcommand{\vsurj}{\downarrow\!\!\!\!\!\makebox[.7pt]{}
                    \raisebox{.3ex}{$\downarrow$}}
\newcommand{\keps}{\varepsilon}
\newcommand{\ku}{\underline}
\newcommand{\kgeq}{\succeq}
\newcommand{\kI}{{\cal J}}
\newcommand{\ktI}{\tilde{\kI}}
\newcommand{\kd}{\displaystyle}

\newcounter{secnum}

\newcommand{\sect}[1]
 {\protect\section{#1}
  \protect\setcounter{secnum}{\value{section}}
  \protect\setcounter{equation}{0}
 \protect\renewcommand{\theequation}{\mbox{\arabic{secnum}.\arabic{equation}}}}

\begin{document}
\title{The versal Deformation of an isolated toric Gorenstein Singularity}

\author{Klaus Altmann\footnotemark[1]\\
%        \small Fachbereich Mathematik der Humboldt-Universit\"at zu Berlin
%        \vspace{-0.7ex}\\ \small Institut f\"ur reine Mathematik, Ziegelstr.
%%13a,
%        D-10099 Berlin, Germany. \vspace{-0.7ex}\\ \small E-mail:
%        altmann@informatik.hu-berlin.de}
	 \small Dept. of Mathematics, M.I.T., Cambridge, MA 02139, U.S.A.
	 \vspace{-0.7ex}\\ \small E-mail: altmann@math.mit.edu}
\footnotetext[1]{Die Arbeit wurde mit einem Stipendium des DAAD unterst\"utzt.}
\date{}
\maketitle

\begin{abstract}
Given a lattice polytope $Q\subseteq \R^n$, we define an affine scheme
$\bar{\CM}$ that
reflects the possibilities of splitting $Q$ into a Minkowski sum.\\
\par
On the other hand, $Q$ induces a toric Gorenstein singularity $Y$, and we
construct a flat family
over $\bar{\CM}$ with $Y$ as special fiber. In case $Y$ has an isolated
singularity only, this
family is versal.
\end{abstract}

\tableofcontents
\par
\vspace{2ex}

%\include{ready}

%{\large

%%%%%%%%%%%%
%
% Abschnitt 1
%
%%%%%%%%%%%%%%
\sect{Introduction}\label{s1}

%%%%%%%%%%
% (1.1)
%%%%%%%%

\neu{11}
The whole deformation theory of an isolated singularity is encoded in its
so-called
versal deformation.\\
For complete intersection singularities this is a family over a
smooth base space - obtained by certain disturbations of the defining
equations.\\
As
soon as we are leaving this class of singularities, the structure of the family
or
even the base space will be more complicated.
It is well known that the base space might consist of several components or
might be non-reduced.\\
In \zitat{9}{2} we will present a (three-dimensional) example
of a singularity
admitting a fat point as base space of its versal deformation.\\
\par

%%%%%%%%%%%%%%
% (1.2)
%%%%%%%%%%

\neu{12}
For two-dimensional cyclic quotient singularities (coinciding with the
two-dimensional affine toric varieties),
the computations of Arndt, Christophersen,
Koll\'{a}r/ Shephard-Barron,
Riemenschneider, and Stevens
provide a description of the versal family - in particular, number and
dimension
of the components of the reduced base (they are smooth) are computed.\\
Christophersen observed that
% , after a finite base change,
the total spaces over these components
are toric varieties again (cf.\ \cite{Ch}).
This suggests the conjecture that the entire deformation theory of affine
toric varieties keeps inside this category. It should be a
challenge to find the versal deformation, its base space, or the total spaces
over the components by purely combinatorial methods.\\
\par

%%%%%%%%%
% (1.3)
%%%%%%%%%%%

\neu{13}
In the present paper we investigate the case of affine, toric
Gorenstein singularities $Y$ given by some lattice polytope $Q$.\\
\par
In \S \ref{s2} and \S \ref{s3} we start with describing an affine scheme
$\bar{\CM}$ which
seems to be interesting independently from the toric or deformation stuff.
It describes the possibilities of spliting
$Q$ into Minkowski
summands. The underlying reduced space is an arrangment of planes corresponding
to
those Minkowski decompositions involving
summands, that are lattice polytopes themselfs.\\
\par
In \S \ref{s5} we
construct a flat family over $\bar{\CM}$
with
the toric Gorenstein singularity $Y$ induced by $Q$ as special fiber.
Computing the Kodaira-Spencer as well as the obstruction map
shows that, in case that
the singularity is isolated, the family is versal (nevertheless trivial
for $\mbox{dim}\,Q\geq 3$).\\
In the general case, the Kodaira-Spencer map is an isomorphism onto the
homogeneous part of $T^1_Y$ with the most interesting multidegree
(cf.\ Theorem \zitat{6}{2}), and the obstruction map is still injective
(cf.\ Theorem \zitat{7}{2}).\\
\par
On the other hand, this family is embedded in a larger (non-flat) family that
equals
a morphism of affine toric varieties: The base space is given by the cone
$C(Q)$ of
Minkowski summands of positive multiples of $Q$, and the total space comes from
the
tautological cone over $C(Q)$ (cf.\ \S \ref{s4} and \S \ref{s5}). In
particular,
for affine, toric, isolated Gorenstein singularities, Christophersen's
observation
(cf.\ \zitat{1}{2})
keeps true (cf.\ \S \ref{s8}).\\
\par
Through the whole paper, an example accompanies the general theory. Further
examples can be
found in \S \ref{s9}.\\
\par

%%%%%%%%%%
% (1.4)
%%%%%%%%%%%%%%%

\neu{14}
{\em Acknowledgements:} I am very grateful to Duco van Straten
and Theo de Jong for constant
encouragement and valuable hints.\\
This paper was written during a one-year-stay at MIT.
I want to thank Richard Stanley and all the other people who made it possible
for me to work at this very interesting and stimulating place.\\
\par

%%%%%%%%%%%%
%
% Abschnitt 2
%
%%%%%%%%%%%%%%
\sect{The Minkowski scheme of a lattice polytope}\label{s2}

%%%%%%%%%%
% (2.1)
%%%%%%%%%%

\neu{21}
Let $Q \subseteq \R^n$ be a lattice polytope, i.e.\ the vertices are contained
in
$\Z^n$. We will always
assume that the edges do not contain any interior lattice points
(cf.\ \zitat{3}{6}), hence, after choosing orientations they are given by
primitive vectors $d^1,\dots, d^N \in \Z^n$.\\
\par
{\bf Definition:}
{\em For every 2-face $\keps <Q$ we define its sign vector
$\ku{\keps}=(\keps_1,\dots, \keps_N) \in \{0,\pm 1\}^N$ by
\[
\keps_i := \left\{
\begin{array}{cl}
\pm 1 & \mbox{if $d^i$ is an edge of $\keps$}\\
0     & \mbox{otherwise}
\end{array} \right.
\]
such that the oriented edges $\keps_i\cdot d^i$ fit to a cycle along the
boundary
of $\keps$.
This determines $\ku{\keps}$ up to sign, and we choose one of both
possibilities.
In particular, $\sum_i \keps_i d^i =0$.}\\
\par

{\bf Example:}
Let us introduce the following example, which will be continued through the
paper:\\
For $Q$ we take the hexagon
\[
Q_6:= \mbox{Conv}\, \{ (0,0), (1,0), (2,1), (2,2), (1,2), (0,1) \} \subseteq
\R^2.
\vspace{2ex}
\]
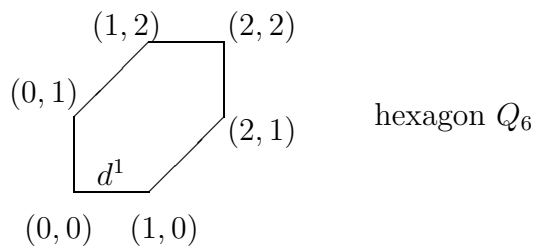
\begin{center}
\unitlength=0.5mm
\linethickness{0.4pt}
\begin{picture}(100.00,64.00)
\put(20.00,20.00){\line(1,0){20.00}}
\put(40.00,20.00){\line(1,1){20.00}}
\put(60.00,40.00){\line(0,1){20.00}}
\put(60.00,60.00){\line(-1,0){20.00}}
\put(40.00,60.00){\line(-1,-1){20.00}}
\put(20.00,40.00){\line(0,-1){20.00}}
\put(16.00,10.00){\makebox(0,0)[cc]{$(0,0)$}}
\put(44.00,10.00){\makebox(0,0)[cc]{$(1,0)$}}
\put(70.00,36.00){\makebox(0,0)[cc]{$(2,1)$}}
\put(70.00,64.00){\makebox(0,0)[cc]{$(2,2)$}}
\put(34.00,64.00){\makebox(0,0)[cc]{$(1,2)$}}
\put(12.00,46.00){\makebox(0,0)[cc]{$(0,1)$}}
\put(30.00,22.00){\makebox(0,0)[cb]{$d^1$}}
\put(100.00,40.00){\makebox(0,0)[lc]{hexagon $Q_6$}}
\end{picture}
\end{center}
Starting with $d^1:= \overline{(0,0) (1,0)}$, the anticlockwise oriented edges
are denoted by
$d^1,\dots, d^6$. As vectors they equal
\[
\begin{array}{ccc}
d^1=(1,0); & d^2=(1,1); & d^3=(0,1);\\
d^4=(-1,0); & d^5=(-1,-1); & d^6=(0,-1).
\end{array}
\]
$Q_6$ is 2-dimensional, hence, it is its own unique 2-face $\keps=Q$. For
$\ku{Q}$ we take
$\ku{Q}=(1,\dots,1)$.\\
\par

%%%%%%%%%%%
% (2.2)
%%%%%%%%%%%

\neu{22}
We define the vector space $V \subseteq \R^N$ by
\[
V:= V(Q):= \{ (t_1,\dots,t_N)\, |\; \sum_i t_i \,\keps_i \,d^i =0\;
\mbox{ for every 2-face } \keps <Q\}.
\]
%({\em Example:} If $n=2$, then $Q$ is a lattice $N$-gon, and $V$ is a
%$(N-2)$-dimensional subspace of $\R^N$.)\\
%\par

Then, $C(Q):= V \cap \R^N_{\geq 0}$ is a rational, polyhedral cone in $V$, and
its
points correspond to the Minkowski summands of positive multiples of $Q$:
Given a point $(t_1,\dots,t_N)\in C(Q)$, the corresponding polytope
$Q_{\ku{t}}$ is built by the
edges $t_i\cdot d^i$ instead of the plain $d^i$ used in $Q$ (cf.\
\zitat{4}{1}).\\
For a Minkowski summand $Q'$ of a positive multiple of $Q$ we will denote its
point
in the cone by $\varrho (Q')\in C(Q)$.\\
\par
({\em Example:} $\varrho(t\cdot Q) = (t,\dots,t) \in C(Q) \subseteq V \subseteq
\R^N.$)\\
\par

%%%%%%%%%%%
% (2.3)
%%%%%%%%%%%

\neu{23}
For each 2-face $\keps <Q$ and for each integer $k\geq 1$
we define the (vector valued)
polynomial
\[
g_{\keps, k}(\underline{t}):= \sum_{i=1}^N t_i^k \, \keps_i \,d^i\,.
\]
Using coordinates of $\R^n$ the $g_{\keps, k}(\underline{t})$ turn into regular
polynomials - for each pair $(\keps,k)$ we will get two linearly independent
ones.\\
\par
We obtain an ideal
\[
\kI := \left( g_{\keps,k} \, | \; \keps < Q, \; k\geq 1 \right)
\subseteq \C [t_1,\dots,t_N],
\]
which defines an affine closed subscheme
\[
\CM := \mbox{Spec}\, ^{\displaystyle \C [\underline{t}]} \!\! \left/ \!\!
_{\displaystyle \kI }\right. \subseteq V_{\C} \subseteq \C^N.
\]
\par

%%%%%%
% Example
%%%%%%%%

{\bf Example:}
For our hexagon $Q_6$ introduced in \zitat{2}{1} we obtain
\[
\kI = \left( t_1^k + t_2^k - t_4^k - t_5^k,\; t_2^k + t_3^k - t_5^k -
t_6^k\,|\; k\geq 1 \right).
\]
\\
Of course, finally many equations are sufficient to generate the ideal $\kI $ -
but we
can even give an effective criterion to see which equations can be dropped:\\
\par
{\bf Proposition:}
{\em Let $\keps <Q$ be a 2-face. Then, $\keps$ is contained in a
two-dimensional
subspace of $\R^n$, and this vector space comes with a natural lattice
(the restriction of
the big lattice $\Z^n$).\\
If $\keps$ is contained in two different strips defined by pairs of
parallel lines of
lattice-distance $\leq k_0$ each, then the equations $g_{\keps,k}\; (k>k_0)$
are
contained
in the ideal generated by $g_{\keps,1}, \dots, g_{\keps,k_0}$.}\\
\par
{\bf Proof:} cf.\ \zitat{3}{3}.\\
\par
{\bf Corollary:}
{\em If $Q$ is contained in $n$ linearly independent strips (defined by pairs
of
parallel hyperplanes) of lattice-thickness $\leq k_0$, then all polynomials
$g_{\keps,k}$ with $k>k_0$ are superfluous.}\\
\par

{\bf Example:}
Obviously, $Q_6$ is contained in at least three strips of thickness 2. Hence,
$\kI $ is
generated in degree $\leq 2$:
\[
\kI = \left( t_1 + t_2 - t_4 - t_5,\quad t_2 + t_3 - t_5 - t_6,\quad
t_1^2 + t_2^2 - t_4^2 - t_5^2,\quad t_2^2 + t_3^2 - t_5^2 - t_6^2 \right) .
\]
\vspace{-1ex}
\\

%%%%%%%%%%%%%
% (2.4)
%%%%%%%%%%

\neu{24}
Denote by $\ell$ the canonical projection
\[
\ell : \C^N \surj \;^{\displaystyle \C^N} \!\!\! \left/
\!_{\displaystyle \C\cdot (1,\dots,1)} \right.
= \;^{\displaystyle \C^N} \!\!\! \left/
\!_{\displaystyle \C\cdot \varrho(Q)} \right. .
\]
On the level of regular functions this corresponds to the inclusion\\
$\C [t_i-t_j\, | \; 1 \leq i,j \leq N ] \subseteq \C[t_1,\dots,t_N]$.\\
\par

{\bf Theorem:}
{\em \begin{itemize}
\item[(1)]
$\kI $ is generated by polynomials from $\,\C [t_i-t_j]$, i.e.\ $\CM =
\ell^{-1}(\bar{\CM})$ for some affine closed subscheme
$\bar{\CM} \subseteq
\,^{\displaystyle V_{\C}} \!\!\! \left/
\!_{\displaystyle \C\cdot \varrho(Q)} \right.
\subseteq
\;^{\displaystyle \C^N} \!\!\! \left/
\!_{\displaystyle \C\cdot \varrho(Q)} \right. .$\\
($\bar{\CM}$ is defined by the ideal $\kI  \cap \C [t_i-t_j]$.)
\item[(2)]
$\kI  \subseteq \C[t_1,\dots,t_N]$ is the smallest ideal (i.e.\ $\CM$ is
the largest closed subscheme of $\C^N$) with
\begin{itemize}
\item[(i)]
property (1) and
\item[(ii)]
containing the ``toric'' equations
%
%\[
%\prod_{i=1}^N t_i^{\displaystyle \langle \keps_id^i, c\rangle_+} -
%\prod_{i=1}^N t_i^{\displaystyle \langle \keps_id^i, c\rangle_-}
%\quad
%\begin{array}[t]{l}
%(\keps < Q \mbox{ is a 2-face of $Q$; $c$ runs through}\\
%\mbox{some $\R$-basis of $\R^n$ contained in }\Z^n).
%\end{array}
%\]
%(For an integer $g$ we denote
%\[
%g_+ := \left\{ \begin{array}{cl}
%g & \mbox{ if } g \geq 0\\
%0 & \mbox{ otherwise}
%\end{array} \right.
%\quad ; \qquad
%g_- := \left\{ \begin{array}{ll}
%0 & \mbox{ if } g \geq 0\\
%-g & \mbox{ otherwise}
%\end{array} \right. .)
%\]
%\end{itemize}
%Moreover, $\kI $ contains all elements $\prod_i t_i^{g^i_+} -
%\prod_i t_i^{g^i_-}$ with\\
%$g \in \Z^N \cap \mbox{span} \left\{
%[\langle\keps _1 d^1,c\rangle, \dots, \langle \keps_N d^N,c \rangle ]\, | \;
%\keps <Q,\, c\in \R^n \right\}$.
%\vspace{2ex}
%\end{itemize}
%
\[
\prod_{i=1}^N t_i^{d_i^+} - \prod_{i=1}^N t_i^{d_i^-}\quad \mbox{ with}
\]
$\ku{d} \in \Z^N \cap \mbox{span} \left\{
[\langle\keps _1 d^1,c\rangle, \dots, \langle \keps_N d^N,c \rangle ]\, | \;
\keps <Q \mbox{ 2-face},\, c\in \R^n \right\}$. \\
(For an integer $h$ we denote
\[
h^+ := \left\{ \begin{array}{cl}
h & \mbox{ if } h \geq 0\\
0 & \mbox{ otherwise}
\end{array} \right.
\quad ; \qquad
h^- := \left\{ \begin{array}{ll}
0 & \mbox{ if } h \geq 0\\
-h & \mbox{ otherwise}
\end{array} \right. .)
\]
\end{itemize}
\end{itemize} }
%\par
{\bf Proof:} cf.\ \zitat{3}{4}.\\
\par

{\bf Example:}
Toric equations for $Q_6$ are for instance
$t_1\,t_2 - t_4\,t_5, \;    t_2\,t_3 - t_5\,t_6,\,$ and $t_1\,t_6 -
t_3\,t_4$.\\
\par

%%%%%%%%%%
% (2.5)
%%%%%%%%%%%%%

\neu{25}
We want to describe the structure of the underlying reduced spaces of $\CM$ or
$\bar{\CM}$.\\
\par
First, we mention the following trivial observations concerning the cone
$C(Q)$:
\begin{itemize}
\item[(i)]
Minkowski summands $Q'$ of $Q$ (instead of a positive multiple of $Q$) are
characterized by the property $\varrho(Q) - \varrho(Q') \in \R^N_{\geq 0}$,
i.e.\ all
components of $\varrho(Q')$ have to be contained in the interval $[0,1]$.
\item[(ii)]
For a Minkowski summand $Q'$ (of some positive multiple of $Q$) the property of
being
a lattice polytope is equivalent to the fact that $\varrho(Q')\in \Z^N$.
\end{itemize}
Now, let $Q=R_0+\dots +R_m$ be a decomposition of $Q$ into a Minkowski sum of
$m+1$
lattice polytopes. Then, the $N$-tuples $\varrho(R_0),\dots,\varrho(R_m)$
consist
of numbers 0 and 1 only, and they sum up to $(1,\dots,1)$. In particular, the
$(m+1)$-plane
$\;\C\cdot\varrho(R_0) + \dots + \C\cdot\varrho(R_m) \subseteq \C^N$
(or its $m$-dimensional image via $\ell$)
is contained in
$\CM$ (in $\bar{\CM}$, respectively).\\
\par

{\bf Remark:}
\begin{itemize}
\item[(1)]
Those $(m+1)$-plane (or its image via $\ell$) is given by the linear equations
$t_i-t_j=0$ (if $d^i, d^j$ belong to a common summand $R_v$).
\item[(2)]
Refinements of Minkowski decompositions (they form a partial ordered set)
correspond to inclusions of the associated planes.
\vspace{2ex}
\end{itemize}
\par

%%%%%%%%%%%
% Theorem
%%%%%%%%%%%

{\bf Theorem:}
{\em $\CM_{\mbox{\footnotesize red}}$ and $\bar{\CM}_{\mbox{\footnotesize
red}}$
equal the union of those flats
corresponding to maximal Minkowski decompositions of $Q$ into lattice
summands.}\\
\par
{\bf Proof:} cf.\ \zitat{3}{5}.\\
\par

{\bf Example:}
$\CM(Q_6)$ and $\bar{\CM}(Q_6)$ are already reduced schemes - for non-reduced
examples
cf.\ \S \ref{s9}. Let us study them directly:
\begin{itemize}
\item
The linear equations allow the following substitution:
\[
\begin{array}{rcl}
t &:=& t_1 \\
s_1 &:=& t_1 - t_3 \\
s_2 &:=& t_4 - t_2 \\
s_3 &:=& t_1 - t_4
\end{array}
\qquad
\begin{array}{rcl}
t_1 &=& t \\
t_2 &=& t-s_2-s_3 \\
t_3 &=& t-s_1 \\
t_4 &=& t-s_3 \\
t_5 &=& t-s_2 \\
t_6 &=& t-s_1-s_3 \, .
\end{array}
\]
\item
The two quadratic equations transform into $s_1\,s_3 = s_2\,s_3 = 0$.
\end{itemize}
In particular, $\bar{\CM}$ is the union of a line and a 2-plane - corresponding
to the Minkowski
decompositions
\[
\begin{array}{rcl}
Q_6 &=& \mbox{Conv}\,\{(0,0), (1,0), (1,1)\} + \mbox{Conv}\,\{(0,0), (0,1),
(1,1)\} \; \mbox{ and}\\
Q_6 &=& \mbox{Conv}\,\{(0,0), (1,0)\} +
\mbox{Conv}\,\{(0,0), (0,1)\} + \mbox{Conv}\,\{(0,0), (1,1)\} .
\end{array}
\vspace{-2ex}
\]
\\
\par

%%%%%%%%%
% (2.6)
%%%%%%%%%

\neu{26}
$\bar{\CM}$ (or $\CM=\ell^{-1}(\bar{\CM})$) reflect the possibilities of
Minkowski decompositions of $Q$:
\begin{itemize}
\item
The underlying reduced space encodes the decompositions of $Q$ into lattice
summands.
\item
Extremal decompositions into rational summands are hidden in the scheme
structure
of $\bar{\CM}$.\\
Its tangent space in 0 (the smallest affine space containing
$\bar{\CM}$) equals $^{\displaystyle V_{\C}}\!\!\!\left/ \!
_{\displaystyle \C \cdot \varrho(Q)} \right.$ - it is the vector space arising
from
the cone $C(Q)$ of Minkowski summands by
killing the summands homothetic to $Q$.
\end{itemize}
Therefore, we will call $\bar{\CM}$ the (affine) {\em Minkowski scheme} of
$Q$.\\
\par

%%%%%%%%%
% Remark
%%%%%%%%%

{\bf Remark:}
The ideals defining $\CM$ and $\bar{\CM}$ are homogeneous. Hence, there are
projective versions of these schemes, too.\\
\par

%%%%%%%%%%%%
%
% Abschnitt 3
%
%%%%%%%%%%%%%%
\sect{Proof of the statements of \S 2}\label{s3}

%%%%%%%%%%
% (3.1)
%%%%%%%%

\neu{31}
Using vectors $c \in \Z^N$ (or selected $c\in \R^N$)
we can evaluate the edges $d^1,\dots,d^N$ to get integers
\[
d_1:= \langle \keps_1 d^1, c \rangle , \dots, d_N:= \langle \keps_N d^N, c
\rangle
\]
for every given 2-face $\keps < Q$.
Doing so, the statements of \S \ref{s2} can be
reduced to much simpler lemmas, which we will present here.\\
\par
Then, all those lemmas are proved using the following recipe:
\begin{itemize}
\item[(i)]
Assume $d_i=\pm 1$ - then the lemmas reduce to well known facts concerning
symmetric functions.
\item[(ii)]
Move to the general case by specialization of variables.
\vspace{2ex}
\end{itemize}
\par

%%%%%%%%%%
% (3.2)
%%%%%%%%

\neu{32}
For the whole \S \ref{s3} we use the following notations:
\vspace{1ex}\\
Let $d_1,\dots,d_N \in \Z$ such that
$\, d_1,\dots, d_M \geq 0, \, d_{M+1},\dots, d_N \leq 0, \,$ and $\sum_{i=1}^N
d_i =0\,$.
%\end{array}$ \\
%$\; \begin{array}[t]{l}
%d_1,\dots, d_M \geq 0, \\
%d_{M+1},\dots, d_N \leq 0, \mbox{ and} \\
%\sum_{i=1}^N d_i =0;
%\end{array}$ \\
\[
\begin{array}{rcccl}
g_k(\underline{t}) &:=& g_{\ku{d},k}(\ku{t}) &:=&
\sum_{i=1}^N d_i \, t_i^k\; , \\
p(\underline{t}) &:=& p_{\ku{d}}(\ku{t}) &:=&
t_1^{d_1}\cdot\dots\cdot t_M^{d_M} - t_{M+1}^{d_{M+1}}\cdot\dots\cdot
t_N^{d_N}.
\end{array}
\]
Denote by $\sigma_k$ and $s_k$ the $k$-th elementary symmetric polynomial and
the sum of the
$k$-th powers of a given set of variables, repectively.
\\
\par

%%%%%%%%%%%%%
% Remark:
%%%%%%%

{\bf Remark:}
For $1\leq i,j \leq M$ or $M+1 \leq i,j \leq N$,
identifying the two variables $t_i$ and $t_j$ (i.e.\ switching from
$\C[\underline{t}]$ to
$^{\displaystyle \C[\underline{t}]} \! /  _{\displaystyle t_i-t_j}$) yields the
following
situation:
\begin{itemize}
\item
$t_i, t_j$ are replaced by a common new variable $\tilde{t}$ (i.e.\ $N$ is
replaced by $N-1$),
\item
$d_i, d_j$ are replaced by $\tilde{d}:= d_i + d_j$, but
\item
$g_k(\underline{t}), p(\underline{t})$ keep their shapes in the new set up.
\end{itemize}
In particular, the general situation can always be obtained via factorization
from the special
case $d_1=\dots = d_M=1;\; d_{M+1}=\dots = d_N=-1$ (and $N=2M$).
Renaming $t_i=x_i, \, t_{M+i}=y_i \,(i\leq M)$ it looks like
\bean
g_k(\underline{x}, \underline{y}) &= &
\left( \sum_{i=1}^M x_i^k \right) - \left( \sum_{i=1}^M y_i^k \right) =
s_k(\underline{x}) -
s_k(\underline{y})\, , \\
p(\underline{x}, \underline{y}) &= &
(x_1\cdot\dots\cdot x_M) - (y_1\cdot\dots\cdot y_M) = \sigma_M(\underline{x}) -
\sigma_M(\underline{y}).
\eean
\\
\par

%%%%%%%%%
% (3.3)
%%%%%%%

\neu{33}
{\bf Lemma:} {\em
If $k_0:= \sum_{i=1}^M d_i = -\sum_{i=M+1}^N d_i$, then the
polynomials $g_k \, (k>k_0)$ are
$\C[\ku{t}]$-linear combinations of the $g_1,\dots,g_{k_0}$. }
(This implies Proposition \zitat{2}{3}.)\\
\par

{\bf Proof:}
As previously discussed, we may regard the special case $d_i=\pm 1$.
In particular, this implies
$k_0=M$.\\
Now, for an arbitrary $k\,(>M)$, the
expression $s_k(\underline{x})$ is a polynomial in either the
$\sigma_1(\underline{x}),\dots, \sigma_M(\underline{x})$ or the
$s_1(\underline{x}), \dots,
s_M(\underline{x})$, say
\[
s_k(\underline{x}) = P_k\left( s_1(\underline{x}),\dots, s_M(\underline{x})
\right) .
\]
Then,
\[
g_k(\underline{x}, \underline{y}) = s_k(\underline{x}) - s_k(\underline{y}) =
P_k\left( s_1(\underline{x}), \dots,s_M(\underline{x})\right) -
P_k\left( s_1(\underline{y}), \dots,s_M(\underline{y})\right),
\]
but for each monomial $s_1^{e_1}\, s_2^{e_2} \dots s_M^{e_M}$ ocuring in $P_k$,
we have
\[
\begin{array}{l}
s_1(\underline{x})^{e_1}\cdot\dots\cdot s_M(\underline{x})^{e_M} -
s_1(\underline{y})^{e_1}\cdot\dots\cdot s_M(\underline{y})^{e_M}  =
\vspace{1ex}\\
\qquad \begin{array}{r}
= \sum_{v=1}^M \sum_{i=1}^{e_v} [s_v(\underline{x}) - s_v(\underline{y})]\cdot
s_1(\underline{x})^{e_1}\dots s_{v-1}(\underline{x})^{e_{v-1}}\,
s_v(\underline{x})^{i-1} \cdot
\qquad\qquad \\
\cdot s_v(\underline{y})^{e_v-i}\, s_{v+1}(\underline{y})^{e_{v+1}} \dots
s_M(\underline{y})^{e_M}
\end{array}
\vspace{1ex}\\
\qquad \begin{array}{r} =
\sum_{v=1}^M g_v(\underline{x}, \underline{y})\cdot \sum_{i=1}^{e_v}
s_1(\underline{x})^{e_1}\dots s_{v-1}(\underline{x})^{e_{v-1}}\,
s_v(\underline{x})^{i-1} \cdot
\qquad\qquad\qquad \\
\cdot s_v(\underline{y})^{e_v-i}\, s_{v+1}(\underline{y})^{e_{v+1}} \dots
s_M(\underline{y})^{e_M},
\end{array}
\end{array}
\]
which proves the lemma.
\hfill$\Box$\\
\par

%%%%%%%%%
%  (3.4)
%%%%%%

\neu{34}
{\bf Lemma:}
{\em \begin{itemize}
\item[(1)]
The ideal $\kI := (g_k\, | \; k\geq 1) \subseteq \C[t_1,\dots, t_N]$ is
generated by polynomials in
$t_i - t_1\; (i=2,\dots,N)$ only.
\item[(2)]
$\kI $ is the smallest ideal generated by polynomials in $t_i-t_1$, which
additionally contains $p$.
\end{itemize} }
(This implies Theorem \zitat{2}{4}.)\\
\par

{\bf Proof:}
(1) Replacing $t_i$ by $t_i-t_1$ as arguments in $g_k$ yields
\bean
g_k(t_1-t_1,\dots,t_N-t_1) &=&
\sum _{i=1}^N d_i\, (t_i-t_1)^k = \sum_{i=1}^N d_i\cdot \left(\sum_{v=0}^k
(-1)^v \, t_1^v \,
t_i^{k-v}\right) \\
&=& \sum_{v=0}^k (-1)^v\, t_1^v \cdot \left( \sum_{i=1}^N d_i\, t_i^{k-v}
\right) =
\sum_{v=0}^k (-1)^v\, t_1^v \, g_{k-v} (\underline{t}).
\eean
In particular, $\left( g_k(\underline{t})\, |\; k\geq 1 \right)$ and
$\left( g_k(\underline{t}- t_1)\, |\; k\geq 1 \right)$ are the same ideals in
$\C[\underline{t}]$.\\
\par
(2) The polynomial rings $\C[\underline{t}]$ and $\C[t_1,\,\underline{t}-t_1]$
are equal, i.e.
each polynomial $q(\underline{t})$ can uniquely be written as
\[
q(\underline{t}) = \sum_{v\geq 0} q_v (t_2-t_1, \dots , t_N-t_1) \cdot t_1^v.
\]
Moreover, if $J\subseteq \C[\underline{t}]$ is an ideal generated by
polynomials in
$\underline{t}-t_1$ only, then for each $q(\underline{t})\in J$ the components
$q_v$ are
automatically contained in $J$, too.\\
\par
Let us determine the components of the polynomial $p$ - we will start with our
special case
again:
\[
p(T+\underline{X}, T+\underline{Y}) = (T+X_1)\cdot \dots \cdot (T+X_M) -
(T+Y_1)\cdot \dots \cdot (T+Y_M)
\]
has $\sigma_k(\underline{X}) - \sigma_k(\underline{Y})$ as coefficient of
$T^{M-k} \;
(k=1,\dots, M)$. Now, there are a polynomial $P_k$ and a non-vanishing rational
number
$c_k$ (not depending on $M$) such that
\[
\sigma_k(\underline{X}) = P_k(s_1(\underline{X}),\dots,s_{k-1}(\underline{X}))
+
c_k\cdot s_k(\underline{X}).
\]
As in the proof of the previous lemma we obtain
\bean
\sigma_k(\underline{X}) - \sigma_k(\underline{Y}) &=&
\begin{array}[t]{r}
P_k(s_1(\underline{X}),\dots,s_{k-1}(\underline{X})) -
P_k(s_1(\underline{Y}),\dots,s_{k-1}(\underline{Y})) + \qquad \\
+ c_k\cdot s_k(\underline{X}) - c_k\cdot s_k(\underline{Y})
\end{array}\\
&=& \sum_{v=1}^{k-1} q_v(\underline{X}, \underline{Y})\cdot g_v(\underline{X},
\underline{Y})
+ c_k \cdot g_k(\underline{X}, \underline{Y})
\eean
for some coefficients $q_v$. Specialization - first by $T\mapsto x_1,\,
X_i\mapsto x_i-x_1,\,
Y_i\mapsto y_i-x_1$, then followed by the usual one - shows that the ideal
generated by the
components $p_v(\underline{t}-t_1)$ of $p$ equals $\kI $.
\hfill$\Box$\\
\par

%%%%%%%%
% (3.5)
%%%%%%

\neu{35}
{\bf Lemma:} {\em
Let $\underline{c}=(c_1,\dots,c_N) \in \C^N$ be a point such that
$g_k(\underline{c})=0$ for each
$k\geq 1$. Then, for every fixed $c\in \C$, we have $\sum_{c_i=c} d_i =0$.}
(This implies Theorem \zitat{2}{5}.)\\
\par

{\bf Proof:}
The equations $\sum_{i=1}^N d_i\, c_i^k =0$ present 0 as a linear combination
of the vectors
$(c_i, c_i^2, c_i^3,\dots)$. On the other hand, it is the Vandermonte that
tells us that this linear
combination has to be a trivial one, i.e.\ the sum of the coefficients $d_i$
belonging to equal
variables vanishes.
\hfill$\Box$\\
\par

%%%%%%%%%
% (3.6)
%%%%%%%

\neu{36}
The polytope $Q$ was assumed to have primitive edges only. Actually, we never
needed this fact
neither in the previous lemmata nor in their proofs. It is only important to
translate these results into
the language of Minkowski summands used in \S \ref{s2}.\\
\par
Droping this condition, similar constructions are possible. However, by
declaring some or all lattice
points contained in edges of $Q$ to be additional, artificial vertices of $Q$,
several possibilities
arise with equal rights. The two extremal cases (add
either no or all possible generalized vertices) seem to
be the most interesting ones.\\
\par
{\bf Remark:}
\begin{itemize}
\item[(1)]
For a natural number $g\in\N$, the polytopes $Q$ (with some fixed set of
possibly artificial vertices)
and $g\cdot Q$ (with the correponding set of vertices) induce the same
Minkowski scheme
$\bar{\CM}$.
\item[(2)]
Let $Q_1\subseteq Q_2$ be the same polytopes with different sets of generalized
vertices.
Then, $\bar{\CM}_1$ is a closed subscheme of $\bar{\CM}_2$. It is defined by
identifying
the variables associated to those generalized edges of $Q_2$, that are
contained in the same
generalized edge of $Q_1$.
\vspace{1ex}
\end{itemize}
{\bf Conjecture:}
Let $Q$ be a lattice polytope such that each extremal Minkowski summand of $Q$
is a lattice
polytope, too. Then, using all generalized vertices of $Q$, the affine schemes
$\CM$ and
$\bar{\CM}$ are reduced.\\
\par
In particular, if $Q$ is an arbitrary lattice polytope (with primitive edges),
then $\bar{\CM}_Q$ would
be embedded in some reduced $\bar{\CM}_{g\cdot Q}$. The non-reduced structure
of
$\bar{\CM}_Q$ would arise as a germ of components visible in $\bar{\CM}_{g\cdot
Q}$ only.\\
\par

%%%%%%%%%%%
%
%  Abschnitt 4
%
%%%%%%%%%%%%%%
\sect{The tautological cone over $C(Q)$}\label{s4}

%%%%%%%%%%
% (4.1)
%%%%%%%%

\neu{41}
In \zitat{2}{2} we have introduced the cone $C(Q)$ of Minkowski summands of
$\R_{\geq 0} \cdot Q$. For an element $(t_1,\dots,t_N)\in C(Q)$ the
corresponding summand
$Q_{\underline{t}}$ was built by the edges $t_i\cdot d^i\; (i=1,\dots,N)$.
However, defining
$Q_{\underline{t}}$ as a particular polytope inside its translation class
requires a closer look:\\
\par
Assume that $0\in \R^n$ coincides with some vertex of the lattice polytope $Q$.
Then, each
vertex $a$ of $Q$ can be reached from there by some walk along the edges of $Q$
- we
obtain
\[
a= \sum_{i=1}^N \lambda_i\, d^i \; \mbox{ for some }
\underline{\lambda}= (\lambda_1,\dots, \lambda_N)\, , \; \lambda_i \in \Z.
\]
Now, given an element $\underline{t}\in C(Q)$, we can define the corresponding
vertex
$a_{\underline{t}}$ (and finally the polytope $Q_{\underline{t}}$ as the convex
hull of all of them)
by
\[
a_{\underline{t}} := \sum_{i=1}^N t_i\, \lambda_i\, d^i.
\]
(The linear equations defining $V= \mbox{span}\,C(Q)$ ensure that this
definition does not depend
on the particular path from $0$ to $a$ through the 1-skeleton of $Q$.)\\
\par

%%%%%%%%
% (4.2)
%%%%%%%%

\neu{42}
{\bf Definition:} {\em
The tautological cone $\tilde{C}(Q) \subseteq \R^n \times V \subseteq \R^{n+N}$
is defined as
\[
\tilde{C}(Q) := \{ (a,\,\underline{t})\, | \; \underline{t}\in C(Q); \, a\in
Q_{\underline{t}} \} .
\] }
\\
{\bf Remark:}
$\tilde{C}(Q)$ is (as $C(Q)$) a rational, polyhedral cone. It is generated by
the pairs
$(a^i_{\underline{t}^j}, \, \underline{t}^j)$ with
\begin{itemize}
\item[$\bullet$]
$a^i$ is a vertex of $Q$ and
\item[$\bullet$]
$\underline{t}^j$ is a fundamental generator of $C(Q)$.
\end{itemize}
(This follows from the simple rule
$(a_{\underline{t}+\underline{t}'},\,\underline{t}+\underline{t}')
= (a_{\underline{t}},\, \underline{t}) + (a_{\underline{t}'},\, \underline{t}')
$ for a vertex $a \in Q$
and $\underline{t}, \underline{t}' \in C(Q)$.)
\vspace{1ex}\\
\par
Defining $\sigma:= \mbox{Cone}(Q) \subseteq \R^{n+1}$ by puting $Q$ into the
hyperplane $(t=1)$,
we obtain a fiber product diagram of rational polyhedral cones:
\[
\begin{array}{ccc}
[\sigma \subseteq \R^{n+1}] & \stackrel{i}{\hookrightarrow} &
[\tilde{C}(Q) \subseteq \R^n \times V ]
\vspace{0.5ex}\\
\vsurj \mbox{\footnotesize pr}_{n+1} && \vsurj \mbox{\footnotesize pr}_V
\vspace{-0.5ex}\\
\R_{\geq 0} \quad& \stackrel{\cdot \varrho (Q)}{\hookrightarrow} & [C(Q)
\subseteq V]
\end{array}
\]
(The vertical maps are projections onto the $(n+1)$-th and the $V$-component,
respectively.
The inclusion $i$ is given by $(t\cdot a;\,t) \mapsto (t\cdot a;\, t,\dots,
t)$.)\\
\par

%%%%%%%%%
% (4.3)
%%%%%%%

\neu{43}
The three cones $\sigma= \mbox{Cone}(Q) \subseteq \R^{n+1},\, \tilde{C}(Q)
\subseteq
\R^n\times V,\,$ and $C(Q)\subseteq V$ define affine toric varieties called $Y,
X,$ and $S$,
respectively. The corresponding rings of regular functions are
$A(Y)=\C[\sigma^{\veee}\cap\Z^{n+1}],\;
A(X)=\C[\tilde{C}(Q)^{\veee}\cap(\Z^n\times V^\ast_{\Z})],$
and $A(S)=\C[C(Q)^{\veee}\cap V^\ast_{\Z}]$.
These varieties come with the following maps:
\begin{itemize}
\item[(i)]
The diagram of \zitat{4}{2} induces a fiber product diagram
\[
\begin{array}{ccc}
Y & \stackrel{i}{\hookrightarrow} & X \\
\downarrow && \downarrow \!\!\pi \\
\C & \hookrightarrow & S\, .
\end{array}
\]
Both horizontal maps are closed embeddings.
(These claims will be checked in \ \zitat{4}{5} and \zitat{4}{8}(1).)
\item[(ii)]
$C(Q) = V \cap \R^N_{\geq 0}$ is contained in $\R^N_{\geq 0}$, and the
inclusion provides a
morphism $p: S \rightarrow \C^N$ defining functions $t_1,\dots, t_N$ on $S$.
The composition $\C \hookrightarrow S \stackrel{p}{\rightarrow} \C^N$ sends $t$
to
$(t,\dots,t)$.
\vspace{2ex}
\end{itemize}
\par

{\bf Remark:}
$Y$ is the affine toric Gorenstein singularity corresponding to the lattice
polytope $Q$. We will
use the map $\pi: X \rightarrow S$ to construct the versal deformation of
$Y$.\\
\par

%%%%%%%%%
% (4.4)
%%%%%%%

\neu{44}
To study the toric varieties $Y, X,$ and $S$ it is important to understand the
dual cones of
$\sigma, \tilde{C}(Q),$ and $C(Q)$, respectively.
%, and the relations between them.
Let us start with the dual cone of $\sigma$:\\
\par
To each non-trivial $c\in\Z^n$ we associate a vertex $a(c)$ of $Q$ and an
integer $\eta_0(c)$
meeting the properties
\bean
\langle Q,\, -c \rangle &\leq& \eta_0(c) \qquad \mbox{ and}\\
\langle a(c),\, -c \rangle &=& \eta_0(c).
\eean
For $c=0$  we define $a(0):=0 \in \R^n$ and $\eta_0(0):=0 \in \Z$.
\pagebreak[3]
\\
\par

{\bf Remark:}
\begin{itemize}
\item[(1)]
With respect to $Q$, $c\neq 0$ is the inner normal vector of the affine
supporting hyperplane
$[\langle \bullet,-c\rangle = \eta_0(c)]$ through $a(c)$. In particular,
$\eta_0(c)$ is uniquely determined,
while $a(c)$ is not.
\item[(2)]
Since $0\in Q$, the integers $\eta_0(c)$ are non-negative.
\vspace{2ex}
\end{itemize}

The dual cone of $\sigma$ is defined as
\[
\sigma^{\veee}:= \{ r\in \R^{n+1}\,|\; \langle\sigma,r\rangle \geq 0\}.
\]
By the definition of $\eta_0$, we have
\[
\partial\sigma^{\veee} \cap \Z^{n+1} = \{[c,\eta_0(c)]\,|\; c\in \Z^n\}\,.
\]
Moreover, if $c^1,\dots,
c^w \in \Z^n\setminus 0$ are those elements producing irreducible pairs
$[c,\eta_0(c)]$ (i.e.\ not
allowing any non-trivial
lattice decomposition $[c,\eta_0(c)]=[c',\eta_0(c')] + [c'',\eta_0(c'')]$),
then the
elements
\[
[c^1,\eta_0(c^1)],\dots, [c^w,\eta_0(c^w)] , \, [\underline{0},1]
\]
form the minimal generator set for $\sigma^{\veee}\cap\Z^{n+1}$ as a semigroup.
Among them are
all pairs $[c,\eta_0(c)]$ corresponding to facets (i.e.\ top dimensional faces)
of $Q$.\\
\par
We obtain a closed embedding $Y\hookrightarrow \C^{w+1}$. The coordinate
functions of
$\C^{w+1}$ will be denoted by $z_1,\dots,z_w,\,t$ corresponding to
$[c^1,\eta_0(c^1)],\dots, [c^w,\eta_0(c^w)], \, [\underline{0},1]$,
respectively.\\
\par

{\bf Example:}
We continue our example $Q_6$ from \S \ref{s2}. Here, the facets of $Q_6$ equal
its edges
$d^1,\dots,d^6$, and they are sufficient for producing
all irreducible pairs\\
$[c^1,\eta_0(c^1)],\dots, [c^6,\eta_0(c^6)]$. We have
\[
\begin{array}{ccccccccc}
c^1 &=& [0,1], & c^2 &=& [-1,1], & c^3 &=& [-1,0], \\
c^4 &=& [0,-1], & c^5 &=& [1,-1], & c^6 &=& [1,0]\, .
\end{array}
\]
The corresponding vertices are (for instance)
\[
a(c^6) = a(c^1) = (0,0), \quad
a(c^2) = a(c^3) = (2,1), \quad
a(c^4) = a(c^5) = (1,2),
\]
and we obtain
\[
%\begin{array}{rclrclrcl}
\eta_0(c^1) = 0,\;  \eta_0(c^2) = 1,\;  \eta_0(c^3) = 2,\;
\eta_0(c^4) = 2,\;  \eta_0(c^5) = 1,\;  \eta_0(c^6) = 0\, .
%\end{array}
\]
\\
\par

%%%%%%%%%
% (4.5)
%%%%%%

\neu{45}
Thinking of $C(Q)$ as a cone in $\R^N$ instead of $V$ allows dualizing
the equation $C(Q)=\R^N_{\geq 0} \cap V$ to get
$C(Q)^{\veee}= \R^N_{\geq 0} + V^\bot$. Hence, for $C(Q)$ as a cone in $V$
we obtain
\[
C(Q)^{\veee}= \;^{\displaystyle \R^N_{\geq 0} + V^\bot}\!\!\left/
_{\displaystyle V^\bot} \right.
= \mbox{Im}\, [\R^N_{\geq 0} \longrightarrow V^\ast].
\]
{\em As already happend with $\R^n$, we do not use different notations
for $\R^N$ and its dual space. However, writing down vectors we try to use
paranthesis and brackets for primal and dual ones, respectively.}\\
\par
The surjection $\R^N_{\geq 0}\surj C(Q)^{\veee}$ induces a map
$\N^N \longrightarrow C(Q)^{\veee}\cap V^\ast_{\Z}$, which does not
need to be
surjective at all. This leads to the following definition:\\
\par
{\bf Definition:}
{\em On $V^\ast_{\Z}$ we introduce a partial ordering ``$\kgeq$'' by
\[
\ku{\eta}\kgeq \ku{\eta}' \quad \Longleftrightarrow
\quad \ku{\eta}-\ku{\eta}' \in \mbox{Im}\, [\N^N\rightarrow
V_{\Z}^\ast] \subseteq C(Q)^{\veee}\cap V^\ast_{\Z}.
\] }
\\
On the geometric level, the non-saturated semigroup
$\mbox{Im}\, [\N^N\rightarrow
V_{\Z}^\ast] \subseteq C(Q)^{\veee}\cap V^\ast_{\Z}$ corresponds to the
scheme theoretical
image ${\bar{S}}$ of $p:S\rightarrow \C^N$, and $S\rightarrow {\bar{S}}$ is its
normalization (cf.\ \zitat{5}{2}).\\
The equations of ${\bar{S}} \subseteq \C^N$ are collected in the kernel of
\[
\C[t_1,\dots,t_N]=\C[\N^N] \stackrel{\varphi}{\longrightarrow}
\C[ C(Q)^{\veee}\cap V^\ast_{\Z}] \subseteq \C[V^\ast_{\Z}],
\]
and it is easy to see that
\bean
\mbox{Ker}\,\varphi &=& \left( \left.
\prod_{i=1}^N t_i^{d_i^+} - \prod_{i=1}^N t_i^{d_i^-}\, \right| \;
\ku{d}\in \Z^N \cap V^\bot \right)\qquad  \mbox{ with}\\
V^\bot &=& \mbox{span}
\left\{ \left.
[\langle\keps _1 d^1,c\rangle, \dots, \langle \keps_N d^N,c \rangle ]\, \right|
\;
\keps <Q \mbox{ is a 2-face},\, c\in \R^n \right\}.
\eean
{\bf Remark:}
Using our new notations, we can reformulate Theorem \zitat{2}{4} now:\\
$\CM \subseteq \C^N$ is the largest closed subscheme that is contained in
${\bar{S}}$ and, additionally, comes from
$^{\displaystyle \C^N} \!\!\! \left/
\!_{\displaystyle \C\cdot \varrho(Q)} \right.$ via $\ell$.\\
\par

On the other hand, dualizing the embedding $\R_{\geq 0} \hookrightarrow C(Q)$
yields
\[
\begin{array}{ccc}
C(Q)^{\veee}\cap V_{\Z}^\ast & \surj & \N\\
\ku{\eta} & \mapsto & \sum_i \eta_i
\end{array}
\]
at the level of semigroups. This map
is surjective, even after restricting to the subset
$\mbox{Im}\, [\N^N\rightarrow V_{\Z}^\ast] $: All vectors $e_i$ corresponding
to the functions $t_i$
map onto $1\in \N$.\\
\par
Geometrically this means that both maps $\C \rightarrow S$ and $\C \rightarrow
{\bar{S}}$ are closed
embeddings,
and the corresponding ideals are
$\left(x^{\ku{\eta}} - x^{\ku{\eta}'}\, \left|\; \ku{\eta}, \ku{\eta}'\in
C(Q)\cap V^\ast_{\Z}
\mbox{ with } \sum_i \eta_i = \sum_i \eta_i' \right. \right)$ and
$(t_i-t_j\,|\; 1\leq i,j \leq N)$,
respectively.
In particular, we got a first contribution to proof the claims made in
\zitat{4}{3}(i).\\
\par

%%%%%%%%%%
% (4.6)
%%%%%%%%

\neu{46}
In the next two sections we take a closer look at the dualized cone
$\tilde{C}(Q)^{\veee}$.\\
\par
{\bf Definition:} {\em
For $c \in \Z^n$ let $\underline{\lambda}^c = (\lambda_1^c,\dots,
\lambda_N^c)\in \Z^N$
describe some path from $0\in Q$ to $a(c)\in Q$ through the 1-skeleton of $Q$
(cf.\ \zitat{4}{1}).
Then,
\[
\underline{\eta}(c) := \left[ -\lambda^c_1\langle d^1,c\rangle, \dots,
-\lambda^c_N\langle d^N,c\rangle \right] \in \Z^N
\]
defines an element $\underline{\eta}(c) \in V_{\Z}^\ast$ not depending on the
choice of the
particular path $\underline{\lambda}^c$.\\ }
\par
(Let $\underline{\tilde{\lambda}}^c$ be a different path from $0$ to $a(c)$ -
it will differ from
$\underline{\lambda}^c$ by some linear combination
$\sum_{\keps < Q} g_{\keps}\, \ku{\keps}$ ($g_{\keps}\in \Z \mbox{ for 2-faces
} \keps < Q$) only.
In particular,
\[
\tilde{\lambda}_i^c \langle d^i,c\rangle - \lambda^c_i \langle d^i,c \rangle =
\sum_{\keps < Q} g_{\keps} \langle \keps_i\,d^i, \, c\rangle,
\]
and we obtain $\underline{\eta}(c)_{\tilde{\lambda}} -
\underline{\eta}(c)_{\lambda} \in V^\bot$.)\\
\par

{\bf Lemma:}
{\em \begin{itemize}
\item[(i)]
$\underline{\eta}(0) = 0 \in V^\ast_{\Z}$.
\item[(ii)]
For all $c\in \Z^n$ we have $\underline{\eta}(c) \kgeq 0$ (in the sense of
Definition \zitat{4}{5}).
\item[(iii)]
$\underline{\eta}$ is convex: $\sum_vg_v\,\underline{\eta}(c^v) \kgeq
\underline{\eta}
(\sum_v g_v\, c^v)$ for natural numbers $g_v\in\N$.
\item[(iv)]
$\sum_{i=1}^N \eta_i(c) = \eta_0(c)$ for arbitrary $c\in \Z^n$.
\vspace{2ex}
\end{itemize} }
\par

{\bf Proof:}
(ii) $a(c)$ is a vertex of $Q$ providing minimal value of the linear function
$\langle \bullet,c \rangle$.
In particular, we can choose a path $\underline{\lambda}^c$ from $0\in Q$ to
$a(c)$ such that this
function decreases in each step, i.e.\ $\lambda_i^c \langle d^i, c \rangle \leq
0 \;
(i=1,\dots,N)$.\\
\par
(iii) We define the following paths through the 1-skeleton of $Q$:
\begin{itemize}
\item
$\underline{\lambda}:=$ path from $0\in Q$ to $a(\sum_v g_v\, c^v)\in Q$,
\item
$\underline{\mu}^v:=$ path from $a(\sum_v g_v\, c^v)\in Q$ to $a(c^v)\in Q$
such that
$\mu_i^v \langle d^i, c^v \rangle \leq 0$ for each $i=1,\dots,N$.
\end{itemize}
Then, $\underline{\lambda}^v := \underline{\lambda} + \underline{\mu}^v$ is a
path from $0\in Q$
to $a(c^v)$, and for $i=1,\dots,N$ we obtain
\bean
\sum_v g_v \, \eta_i (c^v) - \eta_i \left( \sum_v g_v\, c^v \right) &=&
-\sum_v g_v\, (\lambda_i + \mu^v_i)\, \langle d^i,c^v \rangle
+ \lambda_i \left\langle d^i,\, \sum_v g_v\, c^v \right\rangle\\
&=& -\sum_v g_v\, \mu^v_i \, \langle d^i, c^v \rangle \geq 0\, .
\eean
(iv) By definition of $\underline{\lambda}^c$ we have $\sum_{i=1}^N
\lambda_i^c\, d^i = a(c)$.
In particular,
\[
\sum_{i=1}^N \eta_i(c) = - \sum_{i=1}^N \langle \lambda^c_i\, d^i,\,c\rangle =
- \langle a(c),\, c \rangle = \eta_0(c).
\nopagebreak
\vspace{-2ex}
\]
\nopagebreak
\hfill$\Box$
\pagebreak[3]\\
\par

{\bf Example:}
In our hexagon $Q_6$ we choose the following paths from $(0,0)$ to the vertices
$a(c^1),\dots,a(c^6)$, respectively:
\[
\ku{\lambda}^6 = \ku{\lambda}^1 := \ku{0},\quad
\ku{\lambda}^2 = \ku{\lambda}^3 := [1,1,0,0,0,0],\quad
\ku{\lambda}^4 = \ku{\lambda}^5 := [1,1,1,1,0,0]\, .
\]
They provide
\[
\begin{array}{ccccccccc}
\ku{\eta}(c^1) &=& [0,0,0,0,0,0]\,, & \ku{\eta}(c^2) &=& [1,0,0,0,0,0]\,, &
\ku{\eta}(c^3) &=& [1,1,0,0,0,0]\,, \\
\ku{\eta}(c^4) &=& [0,1,1,0,0,0]\,, & \ku{\eta}(c^5) &=& [-1,0,1,1,0,0]\,, &
\ku{\eta}(c^6) &=& [0,0,0,0,0,0]\,.
\end{array}
\]
Since $[1,0,-1,-1,0,1] = [\langle d^1,[1,-1] \rangle,\dots,\langle
d^6,[1,-1]\rangle] \in V^\bot$,
the vector $\ku{\eta}(c^5)$ can be transformed into $[0,0,0,0,0,1]$.\\
\par

{\bf Remark:} The definitions of $a(c), \eta_0(c),$ and $\underline{\eta}(c)$
also make sense for general $c\in \R^n$. Then, $\eta_0(c)\in \R$ and
$\underline{\eta}(c)\in V^\ast$ do not need to be contained in the lattices
anymore.
The previous lemma will keep valid (even for $g_v\in \R_{\geq 0}$ in (iii)),
if the relation ``$\kgeq 0$'' is replaced by the weaker version
``$\in C(Q)^{\veee}$''.\\
\par

%%%%%%%%%
% (4.7)
%%%%%%%

\neu{47}{\bf Proposition:}
{\em \begin{itemize}
\item[(1)]
$\tilde{C}(Q)^{\veee} =
\left\{ \left. [c,\underline{\eta}]\in \R^n \times V^\ast \, \right| \;
\underline{\eta} - \underline{\eta}(c) \in C(Q)^{\veee} \right\}$
\item[(2)]
In particular, $[c, \underline{\eta}(c)] \in \tilde{C}(Q)^{\veee}$, and
moreover, it is the only preimage
of $[c, \eta_0(c)] \in \sigma^{\veee}$ via the surjection $i^{\veee}:
\tilde{C}(Q)^{\veee} \surj
\sigma^{\veee}$.
\item[(3)]
$[c^1,\underline{\eta}(c^1)],\dots, [c^w,\underline{\eta}(c^w)]$ and
$C(Q)^{\veee}\cap V^\ast_{\Z}$
(embedded as $[0,C(Q)^{\veee}]$) generate the semigroup $\tilde{C}(Q)^{\veee}
\cap \left(\Z^n \times V_{\Z}^\ast \right)$. (For recalling the definition of
the $c^1,\dots,c^w$,
cf.\ \zitat{4}{4}.)
\vspace{2ex}
\end{itemize} }
\par

{\bf Proof:}
(1) Let $[c,\underline{\eta}]\in \R^n \times V^\ast$ be given; if some
representative of
$\underline{\eta}$ in $\R^N$ is needed,
then it will be denoted by the same name. We have
the following equivalences:
\bean
[c,\underline{\eta}] \in \tilde{C}(Q)^{\veee} &\Longleftrightarrow&
\langle (Q_{\underline{t}}, \, \underline{t}),\, [c,\underline{\eta}] \rangle
\geq 0 \quad
\mbox{ for each } \underline{t} \in C(Q)\\
&\Longleftrightarrow&
\langle Q_{\underline{t}},\, c\rangle + \langle
\underline{t},\,\underline{\eta} \rangle \geq 0
\quad \mbox{ for each } \underline{t} \in C(Q)\\
&\Longleftrightarrow&
\langle a(c)_{\underline{t}},\, c\rangle + \langle
\underline{t},\,\underline{\eta} \rangle \geq 0
\quad \mbox{ for each } \underline{t} \in C(Q).
\eean
Using some path $\underline{\lambda}^c$ we obtain:
\bean
[c,\underline{\eta}] \in \tilde{C}(Q)^{\veee} &\Longleftrightarrow&
\sum_{i=1}^N t_i\, \lambda_i^c \langle d^i,\,c\rangle + \langle \ku{t},\,
\ku{\eta} \rangle \geq 0
\quad \mbox{ for each } \underline{t} \in C(Q)\\
&\Longleftrightarrow&
\sum_{i=1}^N t_i\cdot \left(\lambda_i^c\, \langle d^i,\,c\rangle + \eta_i
\right) \geq 0
\quad \mbox{ for each } \underline{t} \in C(Q)\\
&\Longleftrightarrow&
\left[ \lambda_1^c\, \langle d^1,\,c\rangle + \eta_1,\dots,
\lambda_N^c\, \langle d^N,\, c \rangle + \eta_N \right] \in C(Q)^{\veee}.
\eean
(2) By part (1) we know that for a $[c,\ku{\eta}]\in\tilde{C}(Q)^{\veee}$ it is
possible to choose
$\R^N$-representatives for
$\ku{\eta}, \ku{\eta}(c)$ such that $\eta_i\geq \eta_i(c)$ for $i=1,\dots,N$.\\
On the other hand, the two equalities $\sum_i \eta_i(c) = \eta_0(c)$ (cf.\ (iv)
of the previous lemma)
and $\sum_i \eta_i = \eta_0(c)$ (corresponding to the fact
$[c,\ku{\eta}] \mapsto [c,\eta_0(c)]$) imply $\ku{\eta} = \ku{\eta}(c)$ then.\\
\par
(3) Let $[c,\ku{\eta}] \in \tilde{C}(Q)^{\veee}$.
Then, $[c,\eta_0(c)]$ is representable as a non-negative
linear combination $[c,\eta_0(c)]=\sum_{v=1}^w p_v\, [c^v, \eta_0(c^v)]$
($p_v\in \N$ if $c\in \Z^n$).
Since both elements $[c, \ku{\eta}(c)]$ and $\sum_vp_v [c^v,\ku{\eta}(c^v)]$
are
preimages of $[c, \eta_0(c)]$ via $i^v$, they must be equal by (2),
and we obtain
\[
\begin{array}{rcl}
[c,\, \ku{\eta}] &=& [c,\, \ku{\eta}(c)] + [0, \ku{\eta} - \ku{\eta}(c)]
\vspace{0.5ex}\\
&=&
\sum_v p_v\,[c^v,\ku{\eta}(c^v)] + [0, \ku{\eta} - \ku{\eta}(c)]\,.
\end{array}
\vspace{-3ex}
\]
\hfill$\Box$\\
\par

%%%%%%%
% (4.8)
%%%%%%%

\neu{48}
Finally, we will take a short look at the geometrical situation reached at this
point.
\begin{itemize}
\item[(1)]
The linear map
\[
\begin{array}{ccc}
\tilde{C}(Q)^{\veee}\cap \left( \Z^n\times V^\ast_{\Z}\right) & \surj
&\sigma^{\veee}\cap\Z^{n+1}\\
\,[c,\, \ku{\eta}] & \mapsto & [c,\,\sum_i\eta_i]
\end{array}
\]
is surjective ($[c,\ku{\eta}(c)]\mapsto [c,\eta_0(c)];\; [0,e_i]\mapsto
[0,1]$). Since
\[
x^{[c,\ku{\eta}]} - x^{[c,\ku{\eta}']} = x^{[c,\ku{\eta}(c)]}\cdot
(x^{[0,\ku{\eta}-\ku{\eta}(c)]} - x^{[0,\ku{\eta}'-\ku{\eta}(c)]}),
\]
the kernel of the corresponding
homomorphism between the semigroup algebras equals the ideal
\[
\left(x^{[0,\ku{\eta}]}-x^{[0,\ku{\eta}']}\, \left| \; \sum_i\eta_i = \sum_i
\eta_i'\right.\right) .
\]
In particular, the map $Y\hookrightarrow X$ is a closed embedding. Moreover,
comparing with the similar statement concerning $C(Q)^{\veee}$ and $\N$ at the
end of
\zitat{4}{5}, we obtain that the diagram of \zitat{4}{3}(i) is a fiber product
diagram, indeed.
\item[(2)]
The elements $[c^1,\ku{\eta}(c^1)], \dots, [c^w,\ku{\eta}(c^w)] \in
\tilde{C}(Q)$ induce some regular
functions $Z_1,\dots,Z_w$ on $X$. They definine a closed embedding
$X \hookrightarrow \C^w \times S$ lifting the ebedding $Y \hookrightarrow
\C^{w+1}$ of
\zitat{4}{4}. Moreover, for $i=1,\dots,N$, $Z_i$ is the only monomial function
lifting $z_i$ from
$Y$ to $X$.
\end{itemize}
We have obtained the following commutative diagram:
\[
\begin{array}{ccccccc}
Y & \hookrightarrow & \C^w\times\C & = & \C^w\times\C\\
\downarrow & \otimes & \downarrow && \downarrow {\scriptstyle \Delta}\\
X & \hookrightarrow & \C^w\times S & \stackrel{p}{\longrightarrow} & \C^w\times
\C^N\\
&& \downarrow & & \downarrow\\
&& S & \stackrel{p}{\longrightarrow} & \C^N & \stackrel{\ell}{\longrightarrow}
& \C^{N-1}\, .
\end{array}
\]
\\
\par

%%%%%%%%%%%
%
%  Abschnitt 5
%
%%%%%%%%%%%%%%
\sect{A flat family over $\bar{\CM}$}\label{s5}

%%%%%%%%%%
% (5.1)
%%%%%%%%

\neu{51}
{\bf Theorem:} {\em Denote by $\bar{X}$ and $\bar{S}$ the scheme theoretical
images of $X$ and
$S$ in $\C^w\times\C^N$ and $\C^N$, respectively. Then,
\begin{itemize}
\item[(1)]
$X\rightarrow \bar{X}$ and $S\rightarrow\bar{S}$ are the normalization maps.
\item[(2)]
$\pi:X\rightarrow S$ induces a map $\bar{\pi}:\bar{X}\rightarrow \bar{S}$, and
$\pi$ can
be recovered from $\bar{\pi}$ via base change $S\rightarrow \bar{S}$.
\item[(3)]
Restricting to $\CM\subseteq \bar{S}$ and composing with $\ell$ turns
$\bar{\pi}$ into a
family
\[
\bar{X} \times_{\bar{S}}\CM \stackrel{\bar{\pi}}{\longrightarrow} \CM
\stackrel{\ell}{\longrightarrow}
\bar{\CM}\,.
\]
It is flat in $0\in \bar{\CM}\subseteq \C^{N-1}$, and the special fiber equals
$Y$.
\vspace{2ex}
\end{itemize} }
The proof of this theorem will fill \S \ref{s5}.\\
\par

%%%%%%%%
% (5.2)
%%%%%%

\neu{52}
The ring of regular functions $A(\bar{S})$ is given as the image of the map
$\C[t_1,\dots,t_N] \rightarrow A(S)$. Since $\Z^N\surj V_{\Z}^\ast$ is
surjective, the rings
$A(\bar{S})\subseteq A(S) \subseteq \C[V^\ast_{\Z}]$ have the same field of
fractions.\\
On the other hand, while $t$-monomials with negative exponents are involved in
$A(S)$, the
surjectivity of $\R^N_{\geq 0}
\surj C(Q)^{\veee}$ tells us that sufficiently high powers of those monomials
always come from $A(\bar{S})$. In particular, $A(S)$ is normal over
$A(\bar{S})$.\\
\par
$A(\bar{X})$ is given as the image
$A(\bar{X})=\mbox{Im}\, (\C[Z_1,\dots,Z_w,t_1,\dots,t_N]\rightarrow A(X))$.
Since $A(X)$ is generated by $Z_1,\dots,Z_w$ over its subring $A(S)$
(cf.\ Proposition \zitat{4}{7}(3)), the same arguments as for $S$ and $\bar{S}$
apply.
Hence, Part (1) of the previous theorem is proved.\\
\par

%%%%%%%%
% (5.3.)
%%%%%%%

\neu{53}
Recalling that $z_1,\dots,z_w,\,t \in A(Y)$ stand for the monomials with
exponents
$[c^1,\eta_0(c^1)],\dots,[c^w,\eta_0(c^w)], [0,1] \in C(Q)^{\vee}\cap
V^\ast_{\Z}$, respectively,
we obtain the following equations defining $Y \subseteq \C^{w+1}$:
\bean
f_{(a,b,\alpha,\beta)}(\ku{z},t) &:=& t^\alpha \,\prod_{v=1}^w z_v^{a_v} -
t^\beta \, \prod_{v=1}^w z_v^{b_v} \\
&& \hspace{-1cm} \mbox{with }
\begin{array}[t]{l}
a,b\in \N^w: \; \sum_v a_v\, c^v = \sum_v b_v\, c^v \quad \mbox{ and}\\
\alpha, \beta \in \N: \; \sum_v a_v\,\eta_0(c^v) + \alpha = \sum_v
b_v\,\eta_0(c^v) + \beta\, .
\vspace{2ex}
\end{array}
\eean

{\bf Example:}
The singularity $Y_6$ induced by the hexagon $Q_6$ equals the cone over the Del
Pezzo
surface of degree 6 (obtained by blowing up three points of $(\PP^2, {\cal
O}(3))$).
As a closed subset of $\C^7$, it is given by the following 9 equations:
\[
\begin{array}{ccc}
f_{(e_1,e_6+e_2,1,0)} = z_1\,t - z_6\,z_2, \; &
f_{(e_2,e_1+e_3,1,0)} = z_2\,t - z_1\,z_3, \; &
f_{(e_3,e_2+e_4,1,0)} = z_3\,t - z_2\,z_4 , \\
f_{(e_4,e_3+e_5,1,0)} = z_4\,t - z_3\,z_5, \; &
f_{(e_5,e_4+e_6,1,0)} = z_5\,t - z_4\,z_6, \; &
f_{(e_6,e_5+e_1,1,0)} = z_6\,t - z_5\,z_1 , \\
f_{(\ku{0},e_1+e_4,2,0)} = t^2 - z_1\,z_4, \; &
f_{(\ku{0},e_2+e_5,2,0)} = t^2 - z_2\,z_5, \; &
f_{(\ku{0},e_3+e_6,2,0)} = t^2 - z_3\,z_6 \, .
\end{array}
\]
\\
\par

%%%%%%%%
% (5.4.)
%%%%%%%

\neu{54}
Defining $c:=\sum_v a_v\, c^v = \sum_v b_v\, c^v$ we can lift the equations of
$Y$ to the following
elements of $\,\C[Z_1,\dots,Z_w,\,t_1,\dots,t_N] \surj
A(\bar{S})[Z_1,\dots,Z_w]$:
\[
F_{(a,b,\alpha,\beta)}(\ku{Z},\ku{t}) := f_{(a,b,\alpha,\beta)}(\ku{Z},t_1) -
\ku{Z}^{[c,\ku{\eta}(c)]} \cdot
\left( \ku{t}^{\alpha e_1 +\sum_v a_v \ku{\eta}(c^v)} -
\ku{t}^{\beta e_1 +\sum_v b_v \ku{\eta}(c^v)} \right) \cdot
\ku{t}^{-\ku{\eta}(c)}\, .
\]
{\bf Remark:}
\begin{itemize}
\item[(1)]
The symbol $\ku{Z}^{[c,\ku{\eta}(c)]}$ means $\prod_{v=1}^w Z_v^{p_v}$ with
natural numbers
$p_v\in \N$ such that
$[c,\ku{\eta}(c)] = \sum_v p_v\, [c^v, \ku{\eta}(c^v)]$ or equivalently
$[c,\eta_0(c)] = \sum_v p_v\, [c^v, \eta_0(c^v)]$. This condition does not
determine the
coefficients $p_v$ uniquely - choose one of the possibilities.
\item[(2)]
By part (iii) of Lemma \zitat{4}{6}, we have $\sum_v a_v\ku{\eta}(c^v),\,
\sum_v b_v\ku{\eta}(c^v)
\kgeq \ku{\eta}(c)$. In particular, representatives of the $\ku{\eta}$'s can be
chosen such that
all $t$-exponents occuring in monomials of $F$ are non-negative, i.e.
$F \in A(\bar{S})[Z_1,\dots,Z_w]$.
\item[(3)]
Mapping $F$ to $A(X)= \oplus_{[c,\ku{\eta}]} \,\C x^{[c,\ku{\eta}]}$
($[c,\ku{\eta}]$ runs through
all elements of $\tilde{C}(Q)^{\veee} \cap (\Z^n\times V_{\Z}^\ast)$;
$Z_v \mapsto x^{[c^v\!, \ku{\eta}(c^v)]},\, t_i \mapsto x^{[0,e_i]}$) yields
\bean
F_{(a,b,\alpha,\beta)} &=&
\begin{array}[t]{r}
\left( t_1^\alpha \, \prod_v Z_v^{a_v} -
\ku{Z}^{[c,\ku{\eta}(c)]} \, \ku{t}^{\alpha e_1 + \sum_v a_v\ku{\eta}(c^v) -
\ku{\eta}(c)} \right) - \qquad\\
- \left( t_1^\beta \, \prod_v Z_v^{b_v} -
\ku{Z}^{[c,\ku{\eta}(c)]} \, \ku{t}^{\beta e_1 + \sum_v b_v\ku{\eta}(c^v) -
\ku{\eta}(c)} \right)
\end{array}\\
& \mapsto &
\begin{array}[t]{r}
\left( x^{\alpha [0,e_1] + \sum_v a_v [c^v\!,\ku{\eta}(c^v)]} -
x^{[c,\ku{\eta}(c)] + \alpha [0,e_1] + \sum_v a_v [0,\ku{\eta}(c^v)] -
[0,\ku{\eta}(c)]} \right) - \quad\\
- \left( x^{\beta [0,e_1] + \sum_v b_v [c^v\!,\ku{\eta}(c^v)]} -
x^{[c,\ku{\eta}(c)] + \beta [0,e_1] + \sum_v b_v [0,\ku{\eta}(c^v)] -
[0,\ku{\eta}(c)]} \right)
\end{array}\\
& = & \; 0 \; - \; 0\, = \,0 \,,
\eean
i.e.\ the polynomials $F_{(a,b,\alpha,\beta)}$ are equations
for $\bar{X} \subseteq \C^w \times \bar{S}$.
\vspace{2ex}
\end{itemize}
\par

{\bf Example:}
In the hexagon example, we obtain the following liftings:
\[
\begin{array}{rclcl}
F_{(e_1,e_6+e_2,1,0)} &=& (Z_1\,t_1 - Z_6\,Z_2 )- Z_1(t_1-t_1) &=&
Z_1\,t_1 - Z_6\,Z_2, \\
F_{(e_2,e_1+e_3,1,0)} &=&( Z_2\,t_1 - Z_1\,Z_3) -Z_2(t_1^2-t_1\,t_2)\,t_1^{-1}
&=&
Z_2\,t_2 - Z_1\,Z_3, \\
F_{(e_3,e_2+e_4,1,0)} &=& (Z_3\,t_1 - Z_2\,Z_4) - Z_3(t_1^2t_2-t_1\,t_2\,t_3)\,
t_1^{-1}t_2^{-1} &=&
Z_3\,t_3 - Z_2\,Z_4, \\
F_{(e_4,e_3+e_5,1,0)} &=& (Z_4\,t_1 - Z_3\,Z_5 )- Z_4(t_1\,t_2\,t_3 -
t_2\,t_3\,t_4)\, t_2^{-1}t_3^{-1}
&=& Z_4\,t_4 - Z_3\,Z_5, \\
F_{(e_5,e_4+e_6,1,0)} &=& (Z_5\,t_1 - Z_4\,Z_6) -
Z_5(t_1\,t_6-t_2\,t_3)\,t_6^{-1} &=&
Z_5\,t_5 - Z_4\,Z_6 ,\\
%(\mbox{using } t_2\,t_3 = t_5\,t_6 \mbox{ in } A(S)) \\
F_{(e_6,e_5+e_1,1,0)} &=& (Z_6\,t_1 - Z_5\,Z_1) - Z_6(t_1-t_6) &=&
Z_6\,t_6 - Z_5\,Z_1 , \\
F_{(\ku{0},e_1+e_4,2,0)} &=&( t_1^2 - Z_1\,Z_4) - (t_1^2-t_2\,t_3) = t_2\,t_3 -
Z_1\,Z_4 &=&
t_5\,t_6 - Z_1\,Z_4, \\
F_{(\ku{0},e_2+e_5,2,0)} &=& (t_1^2 - Z_2\,Z_5) - (t_1^2 - t_3\,t_4) &=&
t_3\,t_4 - Z_2\,Z_5, \\
F_{(\ku{0},e_3+e_6,2,0)} &=& (t_1^2 - Z_3\,Z_6) - (t_1^2-t_1\,t_2) &=&
t_1\,t_2 - Z_3\,Z_6 \, .
\end{array}
\]
\\
\par

%%%%%%%%
% (5.5)
%%%%%%

\neu{55}
To obtain a complete list of equations defining $\bar{X} \subseteq \C^w \times
\bar{S}$,
we have to regard the kernel of the homomorphism
$A(\bar{S})[Z_1,\dots,Z_w]\surj A(\bar{X}) \subseteq A(X)$. It is generated by
the binomials
\[
\begin{array}[t]{r}
\ku{t}^{\ku{\eta}}\, Z_1^{a_1}\cdot\dots\cdot Z_w^{a_w} -
\ku{t}^{\ku{\mu}}\, Z_1^{b_1}\cdot\dots\cdot Z_w^{b_w} \quad \mbox{ such
that}\hspace{4cm}
\vspace{1ex}\\
\begin{array}[t]{l}
\sum_v a_v [c^v,\ku{\eta}(c^v)] + [0,\ku{\eta}] =
\sum_v b_v [c^v,\ku{\eta}(c^v)] + [0,\ku{\mu}]  \, , \vspace{1ex}\\
\mbox{i.e.\ }
\begin{array}[t]{ll}
\bullet & c:= \sum_v a_v \, c^v = \sum_v b_v\, c^v\\
\bullet & \sum_v a_v\, \ku{\eta}(c^v) + \ku{\eta} = \sum_v b_v\, \ku{\eta}(c^v)
+ \ku{\mu}\, .
\end{array}
\end{array}
\end{array}
\]
However,
\bean
\ku{t}^{\ku{\eta}}\,\ku{Z}^{a} -
\ku{t}^{\ku{\mu}}\,\ku{Z}^{b} &=&
\begin{array}[t]{r}
\ku{t}^{\ku{\eta}}\cdot
\left( \prod_v Z_v^{a_v} - \ku{Z}^{[c ,\ku{\eta}(c)]} \, \ku{t}^{\sum_v a_v
\ku{\eta}(c^v) -
\ku{\eta}(c)} \right) - \qquad\\
- \ku{t}^{\ku{\mu}}\cdot
\left( \prod_v Z_v^{b_v} - \ku{Z}^{[c ,\ku{\eta}(c)]} \, \ku{t}^{\sum_v b_v
\ku{\eta}(c^v) -
\ku{\eta}(c)} \right)
\end{array}\\
& = &
\ku{t}^{\ku{\eta}}\cdot F_{(a,p,0,\alpha)} -
\ku{t}^{\ku{\mu}}\cdot F_{(b,p,0,\beta)}
\eean
with $p\in \N^w$ such that $\sum_v p_v [c^v, \ku{\eta}(c^v)] = [c,
\ku{\eta}(c)]$,
$\alpha = \sum_v a_v \eta_0(c^v) - \eta_0(c)$, and
$\beta = \sum_v b_v \eta_0(c^v) - \eta_0(c)$.\\
In particular, $\mbox{Ker}\,(A(\bar{S})[\ku{Z}] \rightarrow A(X))$ is generated
by the polynomials
$F_{(a,b,\alpha,\beta)}$ introduced in \zitat{5}{4}.\\
\par
{\bf Remark:}
\begin{itemize}
\item[(1)]
The inaccuracy caused by writing $\ku{Z}^{[c,\ku{\eta}(c)]}$ for some
undetermined
$Z_1^{p_1}\cdot\dots\cdot Z_w^{p_w}$ (with $\sum_v p_v [c^v, \ku{\eta}(c^v)] =
[c, \ku{\eta}(c)]$)
does not matter: Choosing other coefficients $q_v$ with the same property
yields
\[
Z_1^{p_1}\cdot\dots\cdot Z_w^{p_w} - Z_1^{q_1}\cdot\dots\cdot Z_w^{q_w} =
F_{(p,q,0,0)}(\ku{Z},\ku{t}) = f_{(p,q,0,0)}(\ku{Z},t)\, .
\]
\item[(2)]
Using exponents $\ku{\eta}, \ku{\mu} \in \Z^N$ (instead of $\N^N$), the
binomials
$\ku{t}^{\ku{\eta}}\,\ku{Z}^{a} -
\ku{t}^{\ku{\mu}}\,\ku{Z}^{b}$ generate the kernel of the map
\[
A(S)[\ku{Z}] = A(\bar{S})[\ku{Z}] \otimes_{A(\bar{S})} A(S) \surj
A(\bar{X})  \otimes_{A(\bar{S})} A(S) \surj A(X)\, .
\]
Since $\ku{Z}^a \otimes \ku{t}^{\ku{\eta}} -
\ku{Z}^b \otimes \ku{t}^{\ku{\mu}} =
\ku{Z}^{[c,\ku{\eta}(c)]} \otimes \left(
\ku{t}^{\sum_v a_v \ku{\eta}(c^v) - \ku{\eta}(c) + \ku{\eta}} -
\ku{t}^{\sum_v b_v \ku{\eta}(c^v) - \ku{\eta}(c) + \ku{\mu}} \right) = 0$
in $A(\bar{X})  \otimes_{A(\bar{S})} A(S)$, this implies that the
surjection $A(\bar{X})  \otimes_{A(\bar{S})} A(S) \surj A(X)$ is
injective, too. In particular, part (2) of our theorem is proved.
\vspace{2ex}
\end{itemize}
\par

%%%%%%%
% (5.6)
%%%%%%%%

\neu{56}
Both the point $0\in \bar{\CM}$ and $\C\subseteq \bar{S}$ are given by
the equations $t_i-t_j=0 \; (1\leq i,j,\leq N)$. Modulo these relations, the
equations $F_{(a,b,\alpha,\beta)}$ of $\bar{X}\subset \C^w\times\bar{S}$
specialize to the equations $f_{(a,b,\alpha,\beta)}$ of $Y\subset
\C^{w+1}$.
In particular, $Y$ is the special fiber of the family $\bar{X}\times_{\bar{S}}
\CM \rightarrow \bar{\CM}$.\\
To show the flatness of this family we have to determine all relations
between the $f_{(a,b,\alpha,\beta)}$'s and lift them to relations between
the $F_{(a,b,\alpha,\beta)}$'s.\\
\par
There are three types of relations between the $f_{(a,b,\alpha,\beta)}$'s:
\begin{itemize}
\item[(i)]
$f_{(a,r,\alpha,\gamma)} + f_{(r,b,\gamma,\beta)} = f_{(a,b,\alpha,\beta)}$\\
with
$\begin{array}[t]{ll}
\bullet & \sum_va_vc^v = \sum_v r_v c^v = \sum b_v c^v\; \mbox{ and}\\
\bullet & \sum_v a_v\eta_0(c^v) +\alpha = \sum_vr_v\eta_0(c^v) +\gamma =
\sum_vb_v\eta_0(c^v) + \beta\, .
\end{array}$\\
For this relation, the same equation between the $F$'s is true.
\item[(ii)]
$t\cdot f_{(a,b,\alpha,\beta)} = f_{(a,b,\alpha+1,\beta+1)}\;$ lifts to
$\;t_1\cdot F_{(a,b,\alpha,\beta)} = F_{(a,b,\alpha+1,\beta+1)}$.
\item[(iii)]
$\ku{z}^r\cdot f_{(a,b,\alpha,\beta)} = f_{(a+r,b+r,\alpha,\beta)}$.
\vspace{1ex}\\
With $c:= \sum_va_vc^v = \sum_vb_vc^v,\; \tilde{c}:= c+ \sum_vr_vc^v$ we
obtain
\[
\begin{array}{l}
\ku{Z}^r\cdot F_{(a,b,\alpha,\beta)} - F_{(a+r,b+r,\alpha,\beta)} =
\vspace{1ex}\\
\qquad\begin{array}[t]{r}
=\ku{Z}^{[\tilde{c},\ku{\eta}(\tilde{c})]}
\cdot
\left( \ku{t}^{ \alpha e_1 + \sum_va_v\ku{\eta}(c^v) +
\sum_vr_v\ku{\eta}(c^v)} - \ku{t}^{ \beta e_1 + \sum_vb_v\ku{\eta}(c^v) +
\sum_vr_v\ku{\eta}(c^v)} \right)
\cdot \ku{t}^{-\ku{\eta}(\tilde{c})} - \\
- \ku{Z}^{[c,\ku{\eta}(c)]}\, \ku{Z}^r\cdot
\left( \ku{t}^{ \alpha e_1 + \sum_va_v\ku{\eta}(c^v)} -
\ku{t}^{ \beta e_1 + \sum_vb_v\ku{\eta}(c^v)} \right)\cdot
\ku{t}^{-\ku{\eta}(c)}
\end{array}
\vspace{1ex}\\
\qquad \begin{array}[t]{r} =
\left( \ku{t}^{ \alpha e_1 + \sum_va_v\ku{\eta}(c^v)-\ku{\eta}(c)} -
\ku{t}^{ \beta e_1 + \sum_vb_v\ku{\eta}(c^v)-\ku{\eta}(c)} \right)\cdot
\hspace{4cm}\\
\left( \ku{t}^{\ku{\eta}(c)+\sum_vr_v\ku{\eta}(c^v) - \ku{\eta}(\tilde{c})}
\ku{Z}^{[\tilde{c},\ku{\eta}(\tilde{c})]} - \ku{Z}^{[c,\ku{\eta}(c)]}
\ku{Z}^r \right)\,.
\end{array}
\end{array}
\]
Now, the inequalities
\[
\sum_va_v\ku{\eta}(c^v),\,\sum_vb_v\ku{\eta}(c^v) \kgeq \ku{\eta}(c)
\;\mbox{ and }\;
\ku{\eta}(c)+\sum_vr_v\ku{\eta}(c^v) - \ku{\eta}(\tilde{c}) \kgeq 0
\]
imply that
\begin{itemize}
\item
the first factor is contained in the ideal defining
$0\in \bar{\CM}$, and
\item
the second factor is an equation of $\bar{X}\subseteq
\C^w\times\bar{S}$ (called $F_{(q,p+r,\xi,0)}$ in \zitat{7}{4}).
\end{itemize}
In particular, we have found a lift for the third relation, too.
\end{itemize}
The proof of Theorem \zitat{5}{1} is complete.\\
\par

%%%%%%%%%
% (5.7)
%%%%%%%%%%

\neu{57}
{\bf Example:}
For $Y_6$, the previously constructed family
is contained in $\C^6\times\C^6\stackrel{\mbox{pr}_2}{\longrightarrow}\,
\C^6/_{\displaystyle \C\cdot (1,\dots,1)}$. Its base space is defined by the 4
equations mentioned at the end of \zitat{2}{3}, and for the total space, the 9
equations of \zitat{5}{4}
have to be added.\\
\par

%%%%%%%%%%%%%%%%
%
%  Abschnitt 6
%
%%%%%%%%%%%%%%
\sect{The Kodaira-Spencer map}\label{s6}

%%%%%%%%%%
% (6.1)
%%%%%%%%

\neu{61}
Denote by $E\subseteq \sigma^{\veee}\cap \Z^{n+1}$ the minimal generating
set
\[
E:= \{[c^1,\eta_0(c^1)],\dots,[c^w,\eta_0(c^w)],[\ku{0},1]\}
\]
mentioned in \zitat{4}{4}.
To each vertex $a^j\in Q$
(or equally named fundamental generator $a^j:=(a^j,1)\in \sigma$)
%; $j=0,\dots,M$)
and each element
$R\in \Z^{n+1}$ we associate the subset
\[
E^R_j:= E_{a^j}^R:= \{r\in E\, |\; \langle a^j,r\rangle < \langle
a^j,R\rangle\} \,.
\vspace{2ex}
\]
{\bf Theorem:} (cf. \cite{T1})
{\em
The vector space $T^1_Y$ of infinitesimal deformations of $Y$ is
$\Z^{n+1}$-graded, and
in degree $-R$ it equals
\bean
T^1_Y(-R) &=&
\left( \left.^{\kd L_{\C}\left( \cup_j
%{j=0} ^M
E_j^R\right)} \! \right/
\!\!  \sum_j
%_{j=0}^M
L_{\C}(E_j^R) \right)^\ast
%\\
%& = &
%\left\{ \begin{array}{ll}
%V_{\C}\,(conv\{a^j|\, \langle a^j,R\rangle =1\})\left/\!
%_{\kd (1,\dots,1)}\right. &
%\mbox{if } \langle a^j, R \rangle \leq 1 \;\forall j=0,\dots,M\\
%0 & \mbox{otherwise }.
%\end{array} \right.
\eean
($L(\dots)$ denotes the vector space of linear relations).
}
\\
\par

%%%%%%%%%%
% (6.2)
%%%%%%%%%%

\neu{62}
There is a special degree $R^\ast=[\ku{0},1]\in \Z^{n+1}$ corresponding to the
affine
hyperplane containing $Q$. The associated subsets of $E$ equal
\bean
E^{R^\ast}_j &=& E\cap (a^j)^\bot\\
&=& \{[c^v,\eta_0(c^v)]\, |\; \langle a^j,-c^v\rangle = \eta_0(c^v)\}\,.
\vspace{1ex}
\eean
In \zitat{4}{6}, for each $c\in \Z^n$, we have defined the linear form
$\ku{\eta}(c)\in V_{\Z}^\ast$. Restricted to the cone $C(Q)$, it
maps $\ku{t}$ to
$\mbox{Max}\langle Q_{\ku{t}},-c\rangle = \langle a(c)_{\ku{t}},-c\rangle$.
This induces the following bilinear map:
\[
\begin{array}{cccccl}
\Phi: & ^{\kd V_{\Z}}\!/\!_{\kd (1,\dots,1)}
& \times & L_{\Z}(E\cap \partial\sigma^{\veee})
& \longrightarrow & \Z\\
&\ku{t} & , & q & \mapsto & \sum_{v,i} t_i\, q_v\, \eta_i(c^v)\,.
\end{array}
\]
(Indeed, for $\ku{t}:=\ku{1}$ we obtain
$\sum_{v,i} q_v\, \eta_i(c^v) = \sum_v q_v\, \eta_0(c^v) = 0$ since $q\in
L_{\Z}(E\cap \partial\sigma^{\veee})$.)\\
Moreover, if $q$ comes from one of the submodules
$L_{\Z}(E_j^{R^\ast})\subseteq L_{\Z}(E\cap\partial\sigma^{\veee})$,
we obtain
\bean
\Phi(\ku{t},q) &=&
\sum_v q_v\cdot \mbox{Max}\langle Q_{\ku{t}},-c^v\rangle =
\sum_v q_v\cdot \langle a_{\ku{t}}^j,-c^v\rangle\\
& = &
\langle a_{\ku{t}}^j,-\sum_v q_vc^v\rangle = 0\,.
\vspace{2ex}
\eean
\par

{\bf Theorem:}
{\em
The Kodaira-Spencer map of the family $\bar{X}\times_{\bar{S}}{\cal M}
\rightarrow
\bar{\CM}$ of \S \ref{s5} equals the map
\vspace{-1ex}
\[
T_0\bar{\CM} = ^{\kd V_{\C}}\!\!/\!_{\kd (1,\dots,1)}\longrightarrow
\left( \left. ^{\kd L_{\C}(E\cap\partial\sigma^{\veee})} \! \right/
\! \sum_j
%{j=0}^M
L_{\C}(E_j^{R^\ast}) \right)^\ast
=T^1_Y(-R^\ast)
\]
induced by the previous pairing.  Moreover, this map is an isomorphism.
}
\\
\par

{\bf Proof:}
Using the same symbol $\kI$ for the ideal $\kI\subseteq \C[t_1,\dots,t_N]$
and the intersection $\kI\cap \C[t_i-t_j\,|\; 1\leq i,j\leq N]$
(cf. \zitat{2}{4}), our family
corresponds to the flat $\C[t_i-t_j]/_{\kd \kI} $-module
$\C[\ku{Z},\ku{t}]/_{\kd (\kI, F_\bullet (\ku{Z}, \ku{t}))}$.\\
\par
Now, we fix a non-trivial tangent vector $\ku{t}^0\in V_{\C}$. Via
$t_i\mapsto t+t_i^0\,\keps$ it induces the infinitesimal family given by
the flat $\C[\keps]/\!_{\kd \keps^2}$-module
\[
A_{\ku{t}^0} := \;\left.^{\kd \C[\ku{z}, t, \keps]}\! \right/ \!
_{\kd (\keps^2, F_\bullet(\ku{z}, t+\ku{t}^0\,\keps))}\,.
\]
To obtain the associated $A(Y)$-linear map $I/_{\kd I^2} \rightarrow A(Y)$
($I:=(f_\bullet (\ku{z},t))$ denotes the ideal of $Y$ in $\C^{w+1}$),
we have to compute the images of $f_\bullet (\ku{z},t)$ in
$\keps \,A(Y) \subseteq
A_{\ku{t}^0}$ and divide them by $\keps$:\\
\par
Using the notations of \zitat{5}{3} and \zitat{5}{4}, in
$A_{\ku{t}^0}$ it holds
\bean
0 &=& F_{(a,b,\alpha,\beta)}(\ku{z}, t+\ku{t}^0\,\keps)\\
&=& \!\begin{array}[t]{l}
f_{(a,b,\alpha,\beta)}(\ku{z}, t+t^0_1\,\keps) -\\
\quad - \ku{z}^{[c,\ku{\eta}(c)]}\cdot
\left( (t+\ku{t}^0\,\keps)^{\alpha e_1 + \sum_v a_v
\ku{\eta}(c^v)-\ku{\eta}(c)}
- (t+\ku{t}^0\,\keps)^{\beta e_1 + \sum_v b_v \ku{\eta}(c^v)-\ku{\eta}(c)}
\right)\,.
\end{array}
\eean
The relation $\keps^2=0$ yields
\[
f_{(a,b,\alpha,\beta)}(\ku{z}, t+t^0_1\keps) =
f_{(a,b,\alpha,\beta)}(\ku{z}, t) + \keps\cdot
(\alpha \,t^{\alpha-1} \,t_1^0 \,\ku{z}^a -
\beta \,t^{\beta-1}\, t_1^0 \,\ku{z}^b)\,,
\]
and similarly we can expand the other terms.
%
%\[
%(t+\ku{t}^0\keps)^{\alpha e_1+\sum_v a_v \ku{\eta}(c^v) - \ku{\eta}(c)}
%=
%t^{\alpha + \sum_v a_v \eta_0(c^v) -\eta_0(c)-1} \cdot
%\left( t + \keps\cdot \left(
%\sum_i t_i^0 (\sum_v a_v\eta_i(c^v) -\eta_i(c)0 + t_1^0\,\alpha
%\right) \right)
%\]
%
Eventually, we obtain
\bean
f_{(a,b,\alpha,\beta)}(\ku{z},t) &=&
\begin{array}[t]{l}
-\keps\,t_1^0\, (\alpha\, t^{\alpha -1}\,\ku{z}^a -
\beta\, t^{\beta -1}\,\ku{z}^b )\,+\,
\keps \, \ku{z}^{[c,\ku{\eta}(c)]}\,
t^{\alpha+\sum_v a_v \eta_0(c^v)-\eta_0(c)-1}\cdot
\vspace{1ex}\\
\qquad\qquad\qquad\qquad\cdot
\left[
t_1^0\, (\alpha-\beta) +
%\vspace{1ex}\\
%+\,t^{\alpha+\sum_v a_v \eta_0(c^v)-\eta_0(c)-1} \cdot
+\sum_i t_i^0\, \left(\sum_v (a_v-b_v)\eta_i(c^v) \right) \right]
\end{array}
\vspace{1ex}\\
&=&
\, \keps\cdot x^{\sum_v a_v [c^v,\eta_0(c^v)] +\alpha-1} \cdot
\left( \sum_i t_i^0\, \left(\sum_v (a_v-b_v)\eta_i(c^v) \right) \right)\,.
\eean
(In $\keps \, A(Y)$ we were able to replace the variables $t$ and $z_i$
by $x^{[\ku{0},1]}$ and $x^{[c^v,\eta_0(c^v)]}$, respectively.)\\
\par
On the other hand, we use Theorem \zitat{3}{4} of \cite{T2}: Fixing
$R^\ast\in \Z^{n+1}$, it is the element of $L_{\C}(E\cap\partial\sigma
^{\veee})^\ast$ given by $q\mapsto \sum_{i,v} t_i^0\,q_v\,\eta_i(c^v)$
that corresponds to the infinitesimal deformation of $T^1_Y(-R^\ast)$
defined by the map
\bean
^{\kd I}\! \left/ _{\kd I^2} \right. & \longrightarrow & A(Y)\\
t^\alpha\,\ku{z}^a - t^\beta\,\ku{z}^b & \mapsto &
\left( \sum_{i,v} t_i^0\,  (a_v-b_v)\eta_i(c^v) \right)
\cdot x^{\sum_v a_v [c^v,\eta_0(c^v)] +\alpha-1}\,.
\vspace{-3ex}
\eean
\hfill$\Box$\\
\par

%%%%%%%%%%%
% (6.3)
%%%%%%%%

\neu{63}
To dicuss the meaning of the homogeneous part $T^1_Y(-R^\ast)$ inside the
whole vector space $T^1_Y$, we have to look at the results of
\cite{Gor}:\\
\par
If $\mbox{dim}\, T^1_Y< \infty$ (for instance, if $Y$ has an isolated
singularity), then
\begin{itemize}
\item[(1)]
$T^1_Y=T^1_Y(-R^\ast)$, but
\item[(2)]
$T^1_Y=0$ for $\mbox{dim}\, Y\geq 4$.
\end{itemize}
In particular, the interesting cases arise from 2-dimensional
lattice polygons $Q$ with primitive edges
only. The corresponding 3-dimensional toric varieties $Y$ have an isolated
singularity, and
the Kodaira-Spencer map $T_0\bar{\CM}\rightarrow T^1_Y$ is an
isomorphism.\\
\par
If $T^1_Y$ has infinite dimension, then this comes from the existence of
infinitly many non-trivial homogeneous pieces $T^1_Y(-R)$. Whenever
$\langle a^j, R\rangle \leq 1$ holds for all vertices $a^j\in Q$, we
have
\[
T^1_Y(-R) = V_{\C}(\mbox{conv}\{a^j\,|\; \langle a^j,R\rangle =1\})\,,
\]
i.e. $T^1_Y(-R)$ equals the vector space of Minkowski summands of some
face of $Q$. ($T^1_Y(-R)=0$ for all other $R\in \Z^{n+1}$.)\\
In particular, $T^1_Y(-R^\ast)$ is a typical, but nevertheless extremal
and perhaps the most interesting part of $T^1_Y$.\\
\par

%%%%%%%%%%
%
%  Abschnitt 7
%
%%%%%%%%%%%%%%
\sect{The obstruction map}\label{s7}

%%%%%%%%%
% (7.1)
%%%%%%%%%

\neu{71}
Dealing with obstructions in the deformation theory of $Y$ involves the
$A(Y)$-module $T^2_Y$. Usually, it is defined in the following way:\\
\par
Let
$
m:= \{ ([a,\alpha],[b,\beta])\in \N^{w+1}\times \N^{w+1}\,|\;
\begin{array}[t]{l}
\sum_v a_v\,c^v=\sum_v b_v\,c^v;\\
\sum_v a_v\,\eta_0(c^v) +\alpha = \sum_v b_v \,\eta_0(c^v) +\beta\}
\end{array}
$\vspace{1ex}\\
denote the set
parametrizing the equations $f_{(a,b,\alpha,\beta)}$ generating the ideal
$I\subseteq \C[\ku{z},t]$ of $Y$.
Then,
\[
{\cal R}:=\mbox{Ker}\left( \C[\ku{z},t]^m \surj I\right)
\]
is the module of linear relations between these equations; it contains the
submodule ${\cal R}_0$ of the so-called Koszul relations.\\
\par
{\bf Definition:}\quad
$T^2_Y:=
\;^{\displaystyle \mbox{Hom}\,(^{\displaystyle {\cal R}}\! / \!
_{\displaystyle {\cal R}_0}, A(Y))} \! \left/ \!
_{\displaystyle \mbox{Hom}\,(\C[\ku{z},t]^m,A(Y))} \right.$\,.\\
\par

Now, we have a similar theorem for $T^2_Y$ as we had in \zitat{6}{1} for
$T^1_Y$;
in particular, we use the
notations introduced there.\\
\par

{\bf Theorem:} (cf. \cite{T2})
{\em
The vector space $T^2_Y$ is $\Z^{n+1}$-graded, and
in degree $-R$ it equals
\[
T^2_Y(-R)=
\left( \frac{\displaystyle
\mbox{Ker}\,\left( \oplus_j L_{\C}(E_j^R) \longrightarrow
L_{\C}(E)\right)}{\displaystyle
\mbox{Im}\, \left( \oplus_{\langle a^i,a^j\rangle<Q} L_{\C}(E_i^R\cap E_j^R)
\rightarrow \oplus_i L_{\C}(E_i^R)\right)} \right)^\ast\,.
\]
}
\par

%%%%%%%%%
% (7.2)
%%%%%%%

\neu{72}
In this section we build up the so-called obstruction map. It
detects all possibile infinitesimal extensions of our family over
$\bar{\CM}$ to a flat family over some larger base space.
We follow the explanation given in \S 4 of \cite{RedFund}.\\
\par
As before,
\[
\kI=(g_{\keps,k}(\ku{t}-t_1)\,|\; \keps<Q,\, k\geq 1) =
(g_{\ku{d},k}(\ku{t}-t_1)\,|\; \ku{d}\in V^\bot\cap\Z^N, \, k\geq 1)
\subseteq \C[t_i-t_j]
\]
denotes the homogeneous ideal
of the base space $\bar{\CM}$. Let
\[
\ktI := (t_i-t_j)_{i,j}\cdot \kI +\kI_1\cdot \C[t_i-t_j]
\subseteq
\C[t_i-t_j\,|\; 1\leq i,j\leq N]\,.
\]
Then, $W:= \,^{\kd \kI}\!\!\! \left/ \! _{\kd \ktI} \right.$ is a finitely
dimensional, $\Z$-graded vector space ($W=\oplus_{k\geq 2} W_k$, and
$W_k$ is generated by the polynomials
$g_{\ku{d},k}(\ku{t}-t_1)$). It comes as the kernel in the exact
sequence
\[
0\rightarrow W \longrightarrow
\,^{\kd \C[t_i-t_j]}\!\!\left/ \! _{\kd \ktI} \right.
\longrightarrow
\,^{\kd \C[t_i-t_j]}\!\!\left/ \! _{\kd \kI} \right.
\rightarrow 0\,.
\]
Identifying $t$ with $t_1$ and $\ku{z}$ with $\ku{Z}$,
the tensor product with $\C[\ku{z},t]$ (over $\C$) yields the important,
exact sequence
\[
0\rightarrow W \otimes_{\C} \C[\ku{z},t]
\longrightarrow
\,^{\kd \C[\ku{Z},\ku{t}]}\!\!\left/ \! _{\kd \ktI\cdot \C[\ku{Z},\ku{t}]}
\right.
\longrightarrow
\,^{\kd \C[\ku{Z},\ku{t}]}\!\!\left/ \! _{\kd \kI\cdot \C[\ku{Z},\ku{t}]}
\right.
\rightarrow 0\,.
\vspace{2ex}
\]
Now, let $s$ be any relation with coefficients in $\C[\ku{z},t]$
between the equations $f_{(a,b,\alpha,\beta)}$,
i.e.
\[
\sum s_{(a,b,\alpha,\beta)} \, f_{(a,b,\alpha,\beta)} =0 \quad
\mbox{in } \C[\ku{z},t]\,.
\]
By flatness of our family (cf. \zitat{5}{6}), the components of $s$
can be lifted to $\C[\ku{Z},\ku{t}]$ obtaining an $\tilde{s}$ such that
\[
\lambda(s):=
\sum \tilde{s}_{(a,b,\alpha,\beta)} \, F_{(a,b,\alpha,\beta)} \mapsto 0 \quad
\mbox{in }
\,^{\kd \C[\ku{Z},\ku{t}]}\!\!\left/ \! _{\kd \kI\cdot \C[\ku{Z},\ku{t}]}
\right.\,.
\]
In particular, each relation $s\in {\cal R}$ induces some element
$\lambda(s)\in
W\otimes_{\C} \C[\ku{z},t]$, which is well defined after the additional
projection to $W\otimes_{\C} A(Y)$.
This procedure describes a certain
element $\lambda\in T^2_Y\otimes_{\C}W= \mbox{Hom}(W^\ast, T^2_Y)$
called the obstruction map.\\
\par

{\bf Theorem:}
{\em
%The family is miniversal, iff the Kodaira-Spencer map is an isomorphism and
%the obstruction map $\lambda: W^\ast\rightarrow T^2_Y$ is injective.
%
The obstruction map $\lambda: W^\ast\rightarrow T^2_Y$ is injective.
}\\
\par

{\bf Corollary:}
{\em
If $\mbox{dim}\, T^1_Y < \infty$,
our family equals the versal deformation of $Y$.
In general, we could say that it is ``versal in degree $-R^\ast$''.
}\\
\par

{\bf Proof:}
In \zitat{6}{2} we have proved that the Kodaira-Spencer map is an isomorphism
(at least onto the homogeneous piece $T^1_Y(-R^\ast)$). By a criterion also
described
in \cite{RedFund}, this fact combined with injectivity of the obstruction map
implies
versality.
\hfill$\Box$\\
\par

The remaining part of \S \ref{s7} contains the proof of the previous theorem.\\
\par

%%%%%%%%%%
% (7.3)
%%%%%%%

\neu{73}
We have to improve the notations of \S \ref{s4} and \S \ref{s5}.
Since $\bar{\CM}
\subseteq \bar{S} \subseteq \C^N$, we were able to use the toric
equations (cf. \zitat{2}{4}) during computations modulo $\kI$. In
particular, the exponents $\ku{\eta}\in \Z^N$ of $\ku{t}$ needed be known
modulo $V^\bot$ only; it was enough to define $\ku{\eta}(c)$ as elements
of $V_{\Z}^\ast$.\\
However, to compute the obstruction map, we have to deal with the smaller
ideal $\ktI\subseteq \kI$. Let us start with refining the definitions
of \zitat{4}{6}:

\begin{itemize}
\item[(i)]
For each vertex $a\in Q$ we
choose the following paths through the
1-skeleton of $Q$:
\begin{itemize}
\item[$\bullet$]
$\ku{\lambda}(a):=$ path from $0\in Q$ to $a\in Q$ .
% $\lambda_i(c) \langle d^i,c \rangle \leq 0$ for each $i=1,\dots,N$.
\item[$\bullet$]
$\ku{\mu}^v(a):=$ path from $a\in Q$ to $a(c^v) \in Q$ such that
$\mu_i^v(a) \langle d^i, c^v \rangle \leq 0$ for each $i=1,\dots,N$.
\item[$\bullet$]
$\ku{\lambda}^v(a):= \ku{\lambda}(a) + \ku{\mu}^v(a)$ is then a path
from $0\in Q$ to $a(c^v)$, which depends on $a$.
\end{itemize}
\item[(ii)]
For each $c\in \Z^n$ we use the vertex $a(c)$ to define
\[
\ku{\eta}^c(c):= \left[ -\lambda_1(a(c))\langle d^1,c\rangle,
\dots, -\lambda_N(a(c)) \langle d^N,c\rangle \right] \in \Z^N
\]
and
\[
\ku{\eta}^c(c^v):= \left[ -\lambda_1^v(a(c))\langle d^1,c^v\rangle,
\dots, -\lambda_N^v(a(c)) \langle d^N,c^v\rangle \right] \in \Z^N\,.
\]
\item[(iii)]
For each $c\in \Z^n$ we fix a representation $c=\sum_v p_v^c\, c^v$
($p_v^c\in \N$) such that $\eta_0(c)=\sum_v p_v^c \,\eta_0(c^v)$. (That
means, $c$ is represented only by those generators $c^v$ that define
faces of $Q$ containing the face defined by $c$ itself.)
\end{itemize}

{\bf Remark:}
Let $a\in \N^w$. Denoting $c:=\sum_v a_vc^v$
we obtain
$\sum_v a_v\,\ku{\eta}^c(c^v) - \ku{\eta}^c(c) \in \N^N$
by arguments as in
Lemma \zitat{4}{6}.\\
Moreover, for the special representation $c=\sum_v p_v^c c^v$,
the equation
$\sum_v p_v^c\,\ku{\eta}^c(c^v) = \ku{\eta}^c(c)$ is true.\\
\par

Now, we improve the definition of the polynomials $F_\bullet(\ku{Z},\ku{t})$
given in \zitat{5}{4}. Let $a,b\in \N^w, \alpha,\beta\in \N$ such that
\[
c:= \sum_v a_v\, c^v = \sum_v b_v\, c^v\quad \mbox{and}\quad
\sum_v a_v\, \eta_0(c^v) + \alpha = \sum_v b_v \, \eta_0(c^v) + \beta\,.
\vspace{-0.5ex}
\]
Then,
\[
F_{(a,b,\alpha,\beta)}(\ku{Z},\ku{t}) := f_{(a,b,\alpha,\beta)}(\ku{Z},t_1) -
\ku{Z}^{\ku{p}^c} \cdot
\left( \ku{t}^{\alpha e_1 +\sum_v a_v \ku{\eta}^c(c^v)-\ku{\eta}^c(c)} -
\ku{t}^{\beta e_1 +\sum_v b_v \ku{\eta}^c(c^v)-\ku{\eta}^c(c)} \right)\, .
\vspace{2ex}
\]
\par

%%%%%%%%%
% (7.4)
%%%%%%%%%

\neu{74}
We have to discuss the same three types of relations as we did in
\zitat{5}{6}. Since there is only one single element $c\in \Z^n$ involved
in the relations (i) and (ii), computing modulo $\ktI$ instead of $\kI$
makes no difference in these cases -
we always obtain $\lambda(s)=0$.\\
\par
Let us regard the relation $s:=\left[\ku{z}^r\cdot f_{(a,b,\alpha,\beta)} -
f_{(a+r,b+r,\alpha,\beta)} =0\right]$ $\;(r\in \N^w$).
We will use the following notations:
\begin{itemize}
\item
$c:=\sum_v a_v \, c^v = \sum_v b_v \, c^v;\quad
\ku{p}:=\ku{p}^c;\quad
\ku{\eta}:= \ku{\eta}^c;$
\item
$\tilde{c}:=\sum_v (a_v+r_v) \, c^v = \sum_v (b_v+r_v) \, c^v =
\sum_v (p_v+r_v) \, c^v; \quad
\ku{q}:=\ku{p}^{\tilde{c}};\quad
\tilde{\ku{\eta}}:= \ku{\eta}^{\tilde{c}};$
\item
$\xi:= \sum_i\left(\left( \sum_v (p_v+r_v)\tilde{\eta}_i(c^v)\right)
-\tilde{\eta}_i(\tilde{c})\right) = \sum_v (p_v+r_v)\eta_0(c^v) -
\eta_0(\tilde{c})\,.$
\end{itemize}
Using the same lifting of $s$ to $\tilde{s}$ as in \zitat{5}{6} yields
\[
\begin{array}{rcl}
\lambda(s)&=&\!\begin{array}[t]{l}
\ku{Z}^r\cdot F_{(a,b,\alpha,\beta)}- F_{(a+r,b+r,\alpha,\beta)} \,-
\vspace{1ex}\\
\qquad\qquad -\,
\left( \ku{t}^{ \alpha e_1 + \sum_va_v\ku{\eta}(c^v)-\ku{\eta}(c)} -
\ku{t}^{ \beta e_1 + \sum_vb_v\ku{\eta}(c^v)-\ku{\eta}(c)} \right)\cdot
F_{(q,p+r,\xi,0)}
\end{array}
\vspace{2ex}\\
&=&\!
\begin{array}[t]{l}
-\ku{Z}^{p+r}\cdot \left( \ku{t}^{\alpha e_1 + \sum_v(a_v-p_v)\ku{\eta}(c^v)}
- \ku{t}^{\beta e_1 + \sum_v(b_v-p_v)\ku{\eta}(c^v)} \right)\, +
\vspace{1ex}\\
+\,
\ku{Z}^q\cdot \left( \ku{t}^{\alpha e_1 +
\sum_v(a_v+r_v-q_v)\ku{\tilde{\eta}}(c^v)}
- \ku{t}^{\beta e_1 + \sum_v(b_v+r_v-q_v)\ku{\tilde{\eta}}(c^v)} \right)\, -
\vspace{1ex}\\
-\,
\left( \ku{t}^{\alpha e_1 + \sum_v(a_v-p_v)\ku{\eta}(c^v)}
- \ku{t}^{\beta e_1 + \sum_v(b_v-p_v)\ku{\eta}(c^v)} \right) \cdot
\left(\ku{Z}^q\, \ku{t}^{\sum_v(p_v+r_v-q_v)\ku{\tilde{\eta}}(c^v)}
-\ku{Z}^{p+r}\right)
\end{array}
\vspace{2ex}\\
&=&\!
\begin{array}[t]{l}
\ku{Z}^q\cdot
\left( \ku{t}^{\alpha e_1 + \sum_v(a_v+r_v-q_v)\ku{\tilde{\eta}}(c^v)}
-\ku{t}^{\alpha e_1 + \sum_v(p_v+r_v-q_v)\ku{\tilde{\eta}}(c^v) +
\sum_v (a_v-p_v)\ku{\eta}(c^v)}
\right) -
\vspace{1ex}\\
\qquad - \ku{Z}^q\cdot\left(
\ku{t}^{\beta e_1 + \sum_v(b_v+r_v-q_v)\ku{\tilde{\eta}}(c^v)}
- \ku{t}^{\beta e_1 + \sum_v(p_v+r_v-q_v)\ku{\tilde{\eta}}(c^v) +
\sum_v (b_v-p_v)\ku{\eta}(c^v)} \right)\,.
\end{array}
\end{array}
\]
As in \zitat{5}{6}(iii), we can see that $\lambda(s)$ vanishes modulo $\kI$ (or
even in
$A(\bar{S})$) - just identify $\ku{\eta}$ and $\ku{\tilde{\eta}}$.\\
% Now, we have to calculate modulo $\tilde{\kI}$.\\
\par

%%%%%%%%%
% (7.5)
%%%%%%%%

\neu{75}
In \zitat{7}{2} we already mentioned the isomorphism
\[
W\otimes_{\C} \C[\ku{z},t] \stackrel{\sim}{\longrightarrow}\,
^{\kd \kI\cdot \C[\ku{Z},\ku{t}]} \!\!\left/ \! _{\kd \ktI\cdot
\C[\ku{Z},\ku{t}]} \right.
\]
obtained by identifying $t$ with $t_1$ and $\ku{z}$ with $\ku{Z}$. Now, with
$\lambda(s)$,
we have obtained an element of the right hand side, which has to be interpreted
as an element of
$W\otimes_{\C} \C[\ku{z},t]$.\\
\par
{\bf Lemma:}
{\em
Let $A,B\in \N^N$ such that $\ku{d}:=A-B\in V^\bot$ (i.e. $\ku{t}^A-\ku{t}^B
\in
\kI\cdot \C[\ku{Z},\ku{t}]$). Then, via the previously mentioned isomorphism,
$\ku{t}^A-\ku{t}^B$ corresponds to the element
\[
\sum_{k\geq 1} c_k\cdot g_{\ku{d},k}(\ku{t}-t_1)\cdot t^{k_0-k}
\in W\otimes_{\C}
\C[\ku{z},t]\,
\]
($k_0:=\sum_iA_i$; $c_k$ are the constants occured in \zitat{3}{4}).
In particular, the coefficients from $W_k$ vanish for $k> k_0$.
}\\
\par
{\bf Proof:}
First, we remark that it is allowed to assume that $A=\ku{d}^+$, $B=\ku{d}^-$,
i.e.
$\ku{t}^A-\ku{t}^B = p_{\ku{d}}(\ku{t})$ (cf. \zitat{3}{2}).
(Otherwise we could write this binomial as
\[
\ku{t}^A - \ku{t}^B = \ku{t}^C\cdot \left( \ku{t}^{\ku{d}^+} -
\ku{t}^{\ku{d}^-}\right)\;(C\in \N^N),
\]
and since
\[
\ku{t}^C = (t_1+[\ku{t}-t_1])^C \equiv t_1^{\sum_i C_i} \pmod{(t_i-t_j)}\,,
\vspace{0.5ex}
\]
we would obtain
\vspace{-0.5ex}
\[
\ku{t}^A-\ku{t}^B \equiv t_1^{\sum_iC_i}\cdot \left( \ku{t}^{\ku{d}^+} -
\ku{t}^{\ku{d}^-}\right) \pmod{\ktI}\,. )
\vspace{2ex}
\]
%
%On the other hand,
%\[
%t_1^{\sum_iC_i}\cdot \left(
%\sum_{k\geq 1} c_k\cdot g_{\ku{d},k}(\ku{t}-t_1)\cdot t^{\sum_i d^+_i-k}
%%\right) =
%\sum_{k\geq 1} c_k\cdot g_{\ku{d},k}(\ku{t}-t_1)\cdot t^{\sum_i A_i-k}\,.\,)
%\]
In \zitat{3}{4} we have seen that
\[
p_{\ku{d}}(\ku{t}) = \sum_{k=1}^{k_0} t_1^{k_0-k}\cdot
\left( \sum_{v=1}^{k-1} q_{v,k}(\ku{t}-t_1)\cdot g_{\ku{d},v}(\ku{t}-t_1) +
c_k\cdot g_{\ku{d},k}(\ku{t}-t_1)\right)
\]
(with $k_0:= \sum_id_i^+$). Since $q_{v,k}(\ku{t}-t_1)\in (t_i-t_j)\cdot
\C[t_i-t_j]$, this implies
\[
p_{\ku{d}}(\ku{t}) \equiv \sum_{k=1}^{k_0} t_1^{k_0-k} \cdot c_k\cdot
g_{\ku{d},k}(\ku{t}-t_1) \pmod{\ktI}\,.
%\vspace{1ex}
\]
On the other hand, for $k>k_0$, Lemma \zitat{3}{3} tells us that
$g_{\ku{d},k}(\ku{t}-t_1)$
is a $\C[t_i-t_j]$-linear combination of the elements
$g_{\ku{d},1}(\ku{t}-t_1),\dots, g_{\ku{d},k_0}(\ku{t}-t_1)$.
Then, the degree $k$ part of the corresponding equation shows
$g_{\ku{d},k}(\ku{t}-t_1)\in \ktI$.
\hfill$\Box$\\
\par
{\bf Corollary:}
{\em
Transfered to $W\otimes_{\C}\C[\ku{z},t]$, the element $\lambda(s)$ equals
\[
\sum_{k\geq 1} c_k \cdot g_{\ku{d},k}(\ku{t}-t_1)\cdot \ku{z}^q \cdot
t^{k_0-k}\quad
\mbox{ with} \begin{array}[t]{rcl}
\ku{d} &:=& \sum_v (a_v-b_v)\cdot \left(\ku{\tilde{\eta}}(c^v)-\ku{\eta}(c^v)
\right)\,,\\
k_0 &:=& \alpha + \sum_v (a_v+r_v)\, \eta_0(c^v) - \eta_0(\tilde{c})\,.
\end{array}
\]
The coefficients vanish for $k>k_0$.
}
\\
\par
{\bf Proof:}
We apply the previous lemma to both summands of the $\lambda(s)$-formula of
\zitat{7}{4}. For the first one we obtain
\bean
\ku{d}^a &=& \!
\begin{array}[t]{l}
[ \alpha e_1 + \sum_v(a_v +r_v -q_v)\, \tilde{\ku{\eta}}(c^v)] \, -
\vspace{0.5ex}\\
\qquad\qquad - \,[ \alpha e_1 + \sum_v(p_v+r_v-q_v) \, \tilde{\ku{\eta}}(c^v)
+\sum_v (a_v-p_v)\, \ku{\eta}(c^v) ]
\end{array}\\
&=&
\sum_v (a_v-p_v)\cdot \left(\ku{\tilde{\eta}}(c^v)-\ku{\eta}(c^v) \right)
\qquad \mbox{and}
\vspace{2ex}\\
k_0 &=&
\sum_i \left( \alpha e_1 +
\sum_v(a_v+r_v-q_v)\,\ku{\tilde{\eta}}(c^v)\right)_i\\
&=& \alpha + \sum_v(a_v+r_v-q_v) \,\eta_0(c^v)\,
=\, \alpha + \sum_v (a_v+r_v)\, \eta_0(c^v) - \eta_0(\tilde{c})\,.
\eean
$k_0$ has the same value for both the $a$- and $b$-summand,
and
\[
\begin{array}{rcl}
\ku{d}=\ku{d}^a-\ku{d}^b &=&
\sum_v (a_v-p_v)\cdot \left(\ku{\tilde{\eta}}(c^v)-\ku{\eta}(c^v) \right) -
\sum_v (b_v-p_v)\cdot \left(\ku{\tilde{\eta}}(c^v)-\ku{\eta}(c^v) \right) \\
&=&
\sum_v (a_v-b_v)\cdot \left(\ku{\tilde{\eta}}(c^v)-\ku{\eta}(c^v) \right)\,.
\end{array}
\vspace{-3ex}
\]
\hfill$\Box$\\
\par

%%%%%%%%%%%%
% (7.6)
%%%%%%%%%%%

\neu{76}
Now, we try to approach the obstruction map $\lambda$ from the opposite
direction.
Using the description of $T^2_Y$ given in \zitat{7}{1} we construct an element
of $T^2_Y\otimes_{\C}W$, that afterwards turns out to equal $\lambda$.\\
\par
For a path $\varrho\in \Z^N$ along the edges of $Q$, we will denote by
\[
\ku{d}(\varrho,c):= [\langle \varrho_1\,d^1,\,c\rangle, \dots,
\langle \varrho_N\,d^N,\,c\rangle]\in \Z^N
\]
the vector showing the behaviour of $c\in\Z^n$ passing each particular edge.
If, moreover, $\varrho$ comes from a closed path, $\ku{d}(\varrho,c)$ is also
contained in
$V^\bot$.\\
On the other hand, for each $k\geq 1$, we can use the $\ku{d}$'s from $V^\bot$
to get
elements $g_{\ku{d},k}(\ku{t}-t_1)\in W_k$ generating this vector space.
Composing both procedures we obtain, for each closed path $\varrho\in \Z^N$, a
map
\[
\begin{array}{cccccl}
g^{(k)}(\varrho,\bullet):& \R^n&\longrightarrow &V^\bot &\longrightarrow &W_k\\
&c&&\mapsto && g_{\ku{d}(\varrho,c), k}(\ku{t}-t_1)\,.
\end{array}
\]
\par

{\bf Remark:}
\begin{itemize}
\item[(1)]
Taking the sum over all 2-faces we get a surjective map
\[
\sum_{\keps<Q} g^{(k)}(\ku{\keps},\bullet): \oplus_{\keps<Q} \C^n \surj W_k\,.
\vspace{-2ex}
\]
\item[(2)]
Let $c\in\Z^n$ (having integer coordinates is very important here).
If $\varrho^1, \varrho^2\in\Z^N$ are two paths each connecting vertices
$a,b\in Q$ such that
\begin{itemize}
\item[$\bullet$]
$|\langle a,c\rangle - \langle b,c \rangle | \leq k-1\;$ and
\item[$\bullet$]
$c$ is monoton along both pathes (i.e.
$\langle \varrho^1_i \,d^i, c \rangle;\;
\langle \varrho^2_i \,d^i, c \rangle
\geq 0$
for $i=1,\dots,N$),
\end{itemize}
then $\varrho^1-\varrho^2\in \Z^N$ will be a closed path yielding
$g^{(k)}(\varrho^1-\varrho^2, \,c)=0$ in $W_k$.
\vspace{2ex}
\end{itemize}
\par

{\bf Proof:}
The reason for (1) is the fact that the elements $\ku{d}(\keps,c)$ ($\keps
<Q\,$ 2-face; $c\in \Z^n$) generate $V^\bot$ as a vector space.\\
For the proof of (2), we look at $\ku{d}:=\ku{d}(\varrho^1-\varrho^2,\,c)$.
Since
$d_i=\langle \varrho^1_i\,d^i,\,c\rangle - \langle \varrho^2_i\,d^i,\,c\rangle$
is the difference of two non-negative integers, we obtain
$d_i^+ \leq \langle \varrho_i^1\,d^i,\, c\rangle$.\\
Hence,
\[
\sum_i d_i^+ \leq \sum_i \langle \varrho^1_i\,d^i,\,c\rangle =
\langle b,c\rangle - \langle a,c\rangle \leq k-1\,,
\]
and as in \zitat{7}{5} we obtain $g_{\ku{d},k}(\ku{t}-t_1)
\in \ktI$ by Lemma \zitat{3}{3}.
\hfill$\Box$
\vspace{2ex}\\
\par

Using the notations introduced in \zitat{6}{1} we obtain for $R:=k\, R^\ast,\;
k\geq 2$
\[
E_j^{kR^\ast}=
\{[c^v,\eta_0(c^v)]\,|\; \langle a^j,c^v\rangle + \eta_0(c^v) \leq k-1\}
\cup \{R^\ast\} \subseteq \sigma^{\veee}\cap\Z^{n+1}\,.
\]
Then, we can define the following linear maps :
\[
\begin{array}{cccl}
\psi_j^{(k)}:& L(E_j^{kR^\ast}) & \longrightarrow & W_k\\
& q & \mapsto & \sum_v q_v\cdot
g^{(k)}\left(\ku{\lambda}(a^j)+\ku{\mu}^v(a^j) -
\ku{\lambda}(a(c^v)),\,c^v\right)\,.
\end{array}
\]
(The $q$-coordinate corresponding to $R^\ast\in E_j^{kR^\ast}$ is not used in
the
definition of $\psi_j^{(k)}$.)
\vspace{1ex}\\
\par

{\bf Lemma:}
{\em
Let $\langle a^i,a^j \rangle <Q$ be an edge of the polyhedron $Q$. Then, on
$L(E_i^{kR^\ast} \cap E_j^{kR^\ast}) = L(E_i^{kR^\ast}) \cap L(E_j^{kR^\ast})$,
the
maps $\psi_i^{(k)}$ and  $\psi_j^{(k)}$ coincide.\\
In particular (cf. Theorem \zitat{7}{1}), the $\psi_j^{(k)}$'s induce a linear
map
$\psi^{(k)}: T^2_Y(-kR^\ast)^{\ast}\rightarrow W_k$.
}\\
\par
{\bf Proof:}
Let $q\in L(E_i^{kR^\ast} \cap E_j^{kR^\ast})$. Moreover, we denote by
$\varrho^{ij}\in\Z^N$ the path consisting of the single edge running from
$a^i$ to $a^j$.\\
Then,
\[
\begin{array}{rcl}
\psi_i^{(k)}(q) - \psi_j^{(k)}(q) &=&
\sum_v q_v\cdot g^{(k)}\left(
\ku{\lambda}(a^i) + \ku{\mu}^v(a^i) -
\ku{\lambda}(a^j) + \ku{\mu}^v(a^j) ,\,c^v\right)
\vspace{1ex}\\
&=&
\! \begin{array}[t]{l}
g^{(k)}\left( \ku{\lambda}(a^i) - \ku{\lambda}(a^j) + \varrho^{ij},\,
\sum_v q_v \, c^v\right) \,+
\vspace{0.5ex}\\
\qquad\qquad +\,
\sum_v q_v\cdot g^{(k)}\left(
\ku{\mu}^v(a^i) - \ku{\mu}^v(a^j) - \varrho^{ij},\,c^v\right)\,,
\end{array}
\end{array}
\]
and both summands vanish for several reasons.
The first one is killed simply by the equality $\sum_v q_v\,c^v=0$. For the
second one we can use (2) of the previous remark:\\
If $q_v\neq 0$, then the assumption about $q$ implies the inequalities
\[
0\leq \langle a^i,c^v\rangle - \langle a(c^v),c^v \rangle\, ; \;
\langle a^j,c^v\rangle - \langle a(c^v),c^v \rangle \leq k-1\,.
\]
Hence, assuming w.l.o.g.
$\langle a^i,c^v\rangle \geq \langle a^j,c^v\rangle$, we
can take
$\varrho^1:= -\ku{\mu}^v(a^j)-\varrho^{ij}$
and
$\varrho^2:= -\ku{\mu}^v(a^i)$
to see that
$g^{(k)}\left(\ku{\mu}^v(a^i) - \ku{\mu}^v(a^j) - \varrho^{ij},
\,c^v\right) =0$.
\hfill$\Box$\\
\par

%%%%%%%%%%
% (7.7)
%%%%%%%

\neu{77}
{\bf Proposition:}
{\em
$\;\sum_{k\geq 1} c_k\,\psi^{(k)}$ equals $\lambda^\ast$,
the adjoint of the obstruction map.
}\\
\par

{\bf Proof:}
In Theorem (3.5) of \cite{T2} we gave a dictionary between the two
$T^2$-formulas mentioned in \zitat{7}{1}. Using this result we can
find an element of $\mbox{Hom}(^{\kd \cal R}\!/\!_{\kd {\cal R}_0},\,
W_k\otimes A(Y))$
representing $\psi^{(k)}\in T^2_Y\otimes W_k$ -
it sends relations of type (i) (cf. \zitat{5}{6}) to 0 and deals with relations
of type (ii) and (iii) in the following way:
\[
[\ku{z}^r\, t^{\gamma}\cdot f_{(a,b,\alpha,\beta)} -
f_{(a+r,b+r,\alpha+\gamma,\beta+\gamma)}=0]  \mapsto
\psi_j^{(k)}(a-b)\cdot x^{\sum_v
(a_v+r_v)[c^v,\eta_0(c^v)]+(\alpha+\gamma-k)R^{\ast}}\,,
\]
if
\[
\langle (Q,1),\, \sum_v (a_v+r_v)\, [c^v,\eta_0(c^v)] + (\alpha + \gamma -k)
R^\ast\rangle \geq 0\,,
\]
and $j$ is such that
\[
\langle (a^j,1),\, \sum_v a_v \,[c^v,\eta_0(c^v)] +
(\alpha -k)R^\ast\rangle <0\,;
\]
otherwise the relation is sent to 0 (in particular, if there is not any $j$
meeting the desired property).
\vspace{0.5ex}\\
\par
On $Q$, the linear forms $c:=\sum_va_v\,c^v$ and
$\tilde{c}=\sum_v(a_v+r_v)c^v$ admit their minimal values at the vertices
$a(c)$
and $a(\tilde{c})$, respectively. Hence,
we can transform the previous formula into
\[
\begin{array}{l}
[\ku{z}^r\, t^{\gamma}\cdot f_{(a,b,\alpha,\beta)} -
f_{(a+r,b+r,\alpha+\gamma,\beta+\gamma)}=0]  \mapsto
\psi_{a(c)}^{(k)}(a-b)\cdot x^{\sum_v (a_v+r_v)[c^v,\eta_0(c^v)]+
(\alpha+\gamma-k)R^{\ast}}
\vspace{1ex}\\
\begin{array}[t]{ll}
\mbox{if } &
\begin{array}[t]{l}
\sum_v(a_v+r_v)\eta_0(c^v)-\eta_0(\tilde{c})+(\alpha+\gamma-k) =\\
\qquad =
\langle (a(\tilde{c}),1),\, \sum_v (a_v+r_v)\, [c^v,\eta_0(c^v)] +
(\alpha + \gamma -k)
R^\ast\rangle \geq 0\,,
\end{array}
\vspace{0.5ex}
\\ &
\begin{array}[t]{l}
\sum_va_v\,\eta_0(c^v)-\eta_0(c)+(\alpha-k) =\\
\qquad =
\langle (a(c),1),\, \sum_v a_v \,[c^v,\eta_0(c^v)] +
(\alpha -k)R^\ast\rangle <0
\end{array}
\end{array}
\end{array}
\]
(and mapping to 0 otherwise).\\
\par
Adding the coboundary $h\in \mbox{Hom}\,(\C[\ku{z},t]^m,\, W_k\otimes A(Y))$
\[
h_{(a,\alpha), (b,\beta)}:=
\left\{ \begin{array}{ll}
\psi^{(k)}_{a(c)}(a-b)\cdot x^{\sum_v a_v [c^v,\eta_0(c^v)] +(\alpha -k)R^\ast}
&
\mbox{for } \sum_v a_v\,\eta_0(c^v)-\eta_0(c)+\alpha\geq k\,,\\
0 & \mbox{otherwise}
\end{array} \right.
\]
does not
change the class in $T^2_Y(-kR^\ast)$ (still representing $\psi^{(k)}$), but
improves the represantative from
$\mbox{Hom}(^{\kd \cal R}\!/\!_{\kd {\cal R}_0},\, W_k\otimes A(Y))$.
It still maps type-(i)-relations to 0, and moreover
\[
\begin{array}{l}
[\ku{z}^r\, t^{\gamma}\cdot f_{(a,b,\alpha,\beta)} -
f_{(a+r,b+r,\alpha+\gamma,\beta+\gamma)}=0]  \mapsto
\vspace{1ex}\\
\quad\mapsto \left\{ \begin{array}{ll}
\left(\psi_{a(c)}^{(k)}(a-b)-
\psi_{a(\tilde{c})}^{(k)}(a-b)\right)
\cdot x^{\sum_v (a_v+r_v)[c^v,\eta_0(c^v)]+
(\alpha+\gamma-k)R^{\ast}} &
\mbox{for }
%\sum_v(a_v+r_v)\eta_0(c^v)+(\alpha+\gamma-k)\geq \eta_0(\tilde{c})\\
k_0+\gamma\geq k\\
0 & \mbox{otherwise}\,
\end{array}
\right.
\end{array}
\]
(with $k_0=\alpha + \sum_v(a_v+r_v)\, \eta_0(c^v)-\eta_0(\tilde{c})$).\\
\par
By definition of $\psi^{(k)}_j$ and $g^{(k)}$ we obtain
\[
\begin{array}{l}
\psi_{a(c)}^{(k)}(a-b)- \psi_{a(\tilde{c})}^{(k)}(a-b)\, =
\vspace{0.5ex}\\
\qquad=\,
\sum_v (a_v-b_v)\cdot g^{(k)}\left(
\ku{\lambda}(a(c)) + \ku{\mu}^v(a(c)) -
\ku{\lambda}(a(\tilde{c})) - \ku{\mu}^v(a(\tilde{c})),\;c^v\right)
\vspace{0.5ex}\\
\qquad=\,
\sum_v (a_v-b_v)\cdot g^{(k)}\left(
\ku{\lambda}^v(a(c)) -
\ku{\lambda}^v(a(\tilde{c})) ,\,c^v\right)
\vspace{0.5ex}\\
\qquad=\,
g_{\ku{d},\,k}(\ku{t}-t_1) \quad
\mbox{ with }
\begin{array}[t]{rcl}
\ku{d} &=&
\sum_v(a_v-b_v)\cdot \ku{d} \left(
\ku{\lambda}^v(a(c)) -
\ku{\lambda}^v(a(\tilde{c})) ,\,c^v\right)\\
&=& \sum_v (a_v-b_v)\cdot \left( \tilde{\ku{\eta}}(c^v)-\ku{\eta}(c^v)
\right)\,,
\end{array}
\end{array}
\]
and this completes our proof. Indeed,
\begin{itemize}
\item
for relations of type (ii)
(i.e. $r=0$; $\gamma=1$) we know $c=\tilde{c}$, hence, those relations map
onto 0;
\item
for relations of type (iii) (i.e. $\gamma=0$) we compare the previous
formula with the result obtained in Corollary \zitat{7}{5} -
the coefficients coincide, and
the monomial
$\ku{z}^q\,t^{k_0-k}\in \C[\ku{z},t]$ maps onto
$x^{\sum_v (a_v+r_v)[c^v,\eta_0(c^v)]+
(\alpha+\gamma-k)R^{\ast}}\in A(Y)$.
\vspace{-3ex}
\end{itemize}
\hfill$\Box$\\
\par

%%%%%%%%%%
% (7.8)
%%%%%%%%

\neu{78}
It remains to show that the summands $\psi^{(k)}$ of $\lambda^\ast$ are
indeed surjective maps from $T^2_Y(-kR^\ast)^\ast$ to $W_k$. We will do
so by composing them with auxiliary surjective maps
$p^k: \oplus_{\keps<Q} \C^n \surj T^2_Y(-kR^\ast)^\ast$
yielding
$\psi^{(k)}\circ p^k = \sum_{\keps<Q} g^{(k)}(\ku{\keps},\bullet)$. Then, the
result
follows from the first part of Remark \zitat{7}{6}.\\
\par
In \S 6 of \cite{T2} we used a short exact sequence of complexes called
\[
0\rightarrow L_{\C}(E^R)_{\bullet} \longrightarrow (\C^{E^R})_{\bullet}
\longrightarrow \mbox{span}_{\C}(E^R)_{\bullet}\rightarrow 0
\]
to obtain from Theorem \zitat{7}{1} an isomorphism
\[
T^2_Y(-R) \cong
\left(
\frac{\kd \mbox{Im}\, [\oplus_{\keps<Q} \C^{n+1}\rightarrow
\oplus_{\langle a^i,a^j\rangle < Q} \C^{n+1}] }
{\kd \mbox{Im}\, [\oplus_{\keps <Q} \mbox{span}_{\C}\,
(\cap_{a^j \in \keps} E_j^R) \rightarrow
\oplus_{\langle a^i,a^j\rangle < Q} \C^{n+1}] }
\right)^\ast\,.
\]
Since $R^\ast=[\ku{0},1]\in E_j^{kR^\ast}$ for $k\geq 2$, the induced
surjective
map $\oplus_{\keps<Q}\C^{n+1}\surj T^2_Y(-kR^\ast)^\ast$ factorizes through
$\oplus_{\keps<Q} \C^{n+1}\!/\!_{\kd\C\cdot R^\ast} =
\oplus_{\keps<Q}\C^n$ yielding the auxiliary map $p^k$ just mentioned.
Taking a closer look at the construction of \cite{T2} \S 6,
we can give an explicit
description of $p^k$; eventually we will be able to compute $\psi^{(k)}
\circ p^k$.\\
\par
Let us fix some 2-face $\keps<Q$. Assume that $d^1,\dots,d^M$ are
its counterclockwise orientated edges, i.e.\ the sign vector $\ku{\keps}$ looks
like $\keps_i=1$ for $i=1,\dots,M$ and $\keps_j=0$ otherwise.
Moreover,
we denote the vertices of $\keps<Q$ by $a^1,\dots,a^M$ such that $d^i$
runs from $a^i$ to $a^{i+1}$ ($M+1:=1$).\\
\par
Starting with a $[c,\eta_0]\in \C^{n+1}$ (and, as just mentioned, only the
$c\in\C^n$ is essential) we have to proceed as follows:
\begin{itemize}
\item[(i)]
For $i=1,\dots,M$ we represent $[c,\eta_0]$ as a linear combination of
elements of $E_i^{kR^\ast}\cap E_{i+1}^{kR^\ast}$. (This corresponds to
the lifting from $\mbox{span}_{\C}(E^R)_{\bullet}$ to
$(\C^{E^R})_{\bullet}$.)
\[
[c,\eta_0] =
\sum_v q_{iv}\, [c^v,\eta_0(c^v)] + q_i\, [\ku{0},1]\,,
\]
and $q_{iv}\neq 0$ implies $[c^v,\eta_0(c^v)]\in E_i^{kR\ast}\cap
E_{i+1}^{kR^\ast}$, i.e.
\[
\langle a^i, c^v\rangle + \eta_0(c^v) \leq k-1\,;\quad
\langle a^{i+1}, c^v\rangle + \eta_0(c^v) \leq k-1\,.
\]
\item[(ii)]
We map the result to $\oplus_{i=1}^M \C^{E_i^{kR^\ast}}$ by taking
succesive differences (corresponding to the application of the
differential in
the complex $(\C^{E^R})_{\bullet}$). The result is automatically
contained in $\mbox{Ker}\left( \oplus_i L(E_i^{kR^\ast})\rightarrow
L(E)\right)$, and its $i$-th summand is the linear relation
\[
\sum_v (q_{i,v}-q_{i-1,v})\cdot [c^v,\eta_0(c^v)] +
(q_i-q_{i-1})\cdot [\ku{0},1]=0\,.
\]
\item[(iii)]
Finally, we apply $\psi^{(k)}$ to obtain
\[
\begin{array}{rcl}
\psi^{(k)}(p^k(c)) &=&
\sum_{i=1}^M \sum_v (q_{i,v}-q_{i-1,v})\cdot
g^{(k)} \left(
\ku{\lambda}(a^i)-\ku{\lambda}(a(c^v)) + \ku{\mu}^v(a^i),\;c^v\right)
\vspace{1ex}\\
&=& \! \begin{array}[t]{l}
\sum_{i,v} g^{(k)} \left(
\ku{\lambda}(a^i)-\ku{\lambda}(a(c^v)) + \ku{\mu}^v(a^i),\;
q_{i,v}\,c^v\right)\,-
\vspace{0.5ex}\\
\qquad -\,
\sum_{i,v} g^{(k)} \left(
\ku{\lambda}(a^{i+1})-\ku{\lambda}(a(c^v)) + \ku{\mu}^v(a^{i+1}),\;
q_{i,v}\,c^v\right)
\vspace{1ex}
\end{array}\\
&=&
\sum_{i,v} g^{(k)} \left(
\ku{\lambda}(a^i) -\ku{\lambda}(a^{i+1})
+ \ku{\mu}^v(a^i)- \ku{\mu}^v(a^{i+1}),\;
q_{i,v}\,c^v\right)\,.
\end{array}
\]
\end{itemize}
Similar to the proof of Lemma \zitat{7}{6} we introduce the path $\varrho^i$
consisting of the single edge $d^i$ only. Then, if $q_{iv}\neq 0$
and w.l.o.g. $\langle a^i,c^v\rangle \geq \langle a^{i+1},c^v\rangle$,
the pair of paths $\ku{\mu}^v(a^i)$ and $\ku{\mu}^v(a^{i+1})+\varrho^i$
meets the assumption of Remark \zitat{7}{6}(2) (cf.\ (i)). Hence,
we can proceed as follows:
\[
\begin{array}{rcl}
\psi^{(k)}(p^k(c)) &=& \!
\begin{array}[t]{l}
\sum_{i,v} g^{(k)} \left(
\ku{\lambda}(a^i) -\ku{\lambda}(a^{i+1}) + \varrho^i,\,q_{iv}\,c^v \right)
\,+
\vspace{0.5ex}\\
\qquad\qquad\qquad +\,\sum_{i,v} g^{(k)} \left(
\ku{\mu}^v(a^i)- \ku{\mu}^v(a^{i+1}) - \varrho^i ,\,q_{iv}\,c^v \right)
\end{array}
\vspace{1ex}\\
&=&
\sum_{i=1}^M g^{(k)}\left(
\ku{\lambda}(a^i) -\ku{\lambda}(a^{i+1}) + \varrho^i,\,
\sum_vq_{iv}\,c^v \right)
\vspace{1ex}\\
&=&
\sum_{i=1}^M g^{(k)}\left(
\ku{\lambda}(a^i) -\ku{\lambda}(a^{i+1}) + \varrho^i,\,
c \right)
\vspace{1ex}\\
&=&
g^{(k)}\left( \sum_{i=1}^M \varrho^i,\,c \right)
\vspace{1ex}\\
&=&
g^{(k)}(\ku{\keps},\,c) \,.
\vspace{0.1ex}
\end{array}
\]
\par
Hence, Theorem \zitat{7}{2} is proven.\\
\par

%%%%%%%%%%%%%%%%
%
%  Abschnitt 8
%
%%%%%%%%%%%%%%
\sect{The components of the reduced versal family}\label{s8}

%%%%%%%%%%
% (8.1)
%%%%%%%%

\neu{81}
The components of the reduced base space $\bar{\CM}_{red}$ correspond to
maximal
decompositions of $Q$
into a Minkowski sum $Q=R_0+\dots+R_m$ with lattice polytopes $R_k\subseteq
\R^n$ as
summands.
Intersections of components are obtained by the finest Minkowski decompositions
of $Q$,
that are coarser than all the involved maximal ones.\\
\par

{\bf Theorem:}
{\em
Fix such a Minkowski decomposition. Then, the corresponding component (or
intersection of
components)
$\bar{\CM}_0$
is isomorphic to
$^{\displaystyle \C^{m+1}}\!\!/\!{\displaystyle \C\cdot (1,\dots,1)}$, and the
restriction $X_0\rightarrow \C^m$
of the versal family can be described as follows:
\begin{itemize}
\item[(i)]
Define the cone
\vspace{-0.5ex}
\[
\tilde{\sigma}:= \mbox{Cone}\,\left( \bigcup_{k=0}^m (R_k\times \{e^k\})
\right)
\subseteq \R^{n+m+1}\, ,
\]
it contains $\sigma = \mbox{Cone}\,(Q\times\{1\})\subseteq \R^{n+1}$ via the
diagonal embedding
$\R^{n+1}\hookrightarrow \R^{n+m+1}\; ((a,1)\mapsto (a;1,\dots,1))$. The
inclusion
$\sigma \subseteq \tilde{\sigma}$ induces a closed embedding of the
affine toric varities defined by these cones - this gives $Y\hookrightarrow
X_0$.
\item[(ii)]
The projection $\R^{n+m+1} \surj \R^{m+1}$ provides $m+1$ regular functions on
$X_0$, i.e.\ we
obtain a map $X_0\rightarrow \C^{m+1}$. Composing this map with\\
$\ell: \C^{m+1}\surj ^{\displaystyle \C^{m+1}}\!\!/\!{\displaystyle \C\cdot
(1,\dots,1)}$ yields the
family.
\vspace{2ex}
\end{itemize} }
% While part (1) was already formulated in Theorem \zitat{2}{5}
% (and proved in \zitat{3}{5}),
We will
use the \zitat{8}{2} and \zitat{8}{3} to prove the theorem.\\
\par

%%%%%%%%%%
% (8.2)
%%%%%%%

\neu{82}
We already know (cf.\ \zitat{2}{5}) that both the space
\[
\bar{\CM}_0 =\; ^{\displaystyle \C^{m+1}}\!\!\!\!/\!_{\displaystyle \C\cdot
(1,\dots,1)}
\;\subseteq \;^{\displaystyle \C^N}\!\!\!\!/\!_{\displaystyle \C\cdot
(1,\dots,1)}
\]
and its pullback $\CM_0\subseteq\C^N$ are
given by the equations $t_i-t_j=0$ (if $d^i,d^j$ belong to a common summand
$R_k$ of $Q$).\\
\par
There is a chain of inclusions
$\CM_0 \subseteq \CM \subseteq \bar{S} \subseteq \C^N$, and
the map $\CM_0 \hookrightarrow \bar{S}$ factorizes through an embedding
$\CM_0 \hookrightarrow S$. It is given by the surjection of $\C$-algebras
\[
\begin{array}{ccc}
\C[C(Q)^{\veee}\cap V_{\Z}^\ast] &\surj& \C[T_0,\dots,T_m]\\
t_i & \mapsto & T_k \mbox{ with } d^i\in R_k
\end{array}
\]
coming from the linear map
\[
\begin{array}{ccccc}
\R^{m+1} & \hookrightarrow & V & \hookrightarrow & \R^N\\
e^k & \mapsto & \sum_{d^i\in R_k} e^i\, .
\end{array}
\]
The matrix of $\R^{m+1}\hookrightarrow \R^N$ equals
the incidence matrix between edges and Minkowski summands of $Q$, and the
space $\CM_0=\C^{m+1}$ corresponds to the cone $\R^{m+1}_{\geq0}=\R^{m+1}\cap
C(Q)$.\\
\par

%%%%%%%%%%%%
% (8.3)
%%%%%%

\neu{83}
The family $X_0\rightarrow\bar{\CM}_0$ arises from $\bar{X}\times_{\bar{S}}\CM
\rightarrow \bar{\CM}$ via base change $\bar{\CM}_0 \hookrightarrow \bar{\CM}$.
We
obtain
\[
\begin{array}{ccccccl}
X_0 &= &
(\bar{X}\times_{\bar{S}}\CM) \times _{\bar{\CM}} \bar{\CM}_0
&=&
\bar{X}\times_{\bar{S}} (\CM \times _{\bar{\CM}} \bar{\CM}_0 )
&=&
\bar{X}\times_{\bar{S}} \CM_0 \\
&=&
\bar{X}\times_{\bar{S}} (S\times_S \CM_0)
&=&
(\bar{X}\times_{\bar{S}} S) \times_S \CM_0
&=&
X\times_S \CM_0\, .
\end{array}
\]
Hence, with $\tilde{\sigma}$ as defined is the theorem, it remains to show that
\[
\begin{array}{ccc}
\C[\tilde{C}(Q)^{\veee} \cap (\Z^n \times V^\ast_{\Z})]
&\stackrel{\psi}{\longrightarrow}&
\C[\tilde{\sigma}^{\veee}\cap \Z^{n+m+1}]\\
\bigcup \!| && \bigcup \!|\\
\C[C(Q)^{\veee}\cap V^\ast_{\Z}] &\surj&
\C[\N^{m+1}] = \C[T_0,\dots,T_m]
\end{array}
\]
is a tensor product diagram:
\begin{itemize}
\item[(i)]
$\tilde{\sigma}$ is the preimage of $\R^{m+1}_{\geq 0}\subseteq C(Q)$ via the
projection
$\tilde{C}(Q)\surj C(Q)$. In particular, $\tilde{\sigma}\subseteq \tilde{C}(Q)$
causes a surjective
map $\psi_{\R}: \tilde{C}(Q)^{\veee} \surj \tilde{\sigma}^{\veee}$.
\item[(ii)]
To show surjectivity at the level of lattices (i.e.\ surjectivity of $\psi$) we
start with an element
$[c,\ku{\eta}]\in \tilde{C}(Q)^{\veee}$ and suppose its image
$\psi_{\R}([c,\ku{\eta}])$ to be contained in
$\tilde{\sigma}^{\veee}\cap \Z^{n+m+1}$.\\
In particular, $c\in \Z^n$, and we obtain
$[c,\ku{\eta}(c)]\in \tilde{C}(Q)^{\veee}\cap(\Z^n\times V_{\Z}^\ast)$
implying that
\[
[0,\ku{\eta}-\ku{\eta}(c)] = [c,\ku{\eta}] - [c,\ku{\eta}(c)] \in
[0,C(Q)^{\veee}] \subseteq
\tilde{C}(Q)^{\veee}
\]
maps to an element of $\N^{m+1}\subseteq \tilde{\sigma}^{\veee}\cap
\Z^{n+m+1}$.\\
On the other hand, surjectivity of $C(Q)^{\veee}\cap V_{\Z}^\ast \surj
\N^{m+1}$ causes that this
element can be reached by some $[0,\ku{\mu}]\in [0,\, C(Q)^{\veee}\cap
V_{\Z}^\ast]$,
too.\\
Hence, $[c,\ku{\eta}(c) + \ku{\mu}] \in
\tilde{C}(Q)^{\veee}\cap(\Z^n\times V_{\Z}^\ast)$
is a lattice-preimage of $\psi_{\R}([c,\ku{\eta}])$.
\item[(iii)]
The same methods applies for showing that $\mbox{Ker}\,\psi$ is generated by
the same elements
as $\mbox{Ker}\,(\C[C(Q)^{\veee}\cap V^\ast_{\Z}] \rightarrow \C[\N^{m+1}])$:\\
If $[c^1,\ku{\eta}^1],\, [c^2,\ku{\eta}^2] \in
\tilde{C}(Q)^{\veee}\cap(\Z^n\times V_{\Z}^\ast)$ have the
same image in $\C[\tilde{\sigma}^{\veee}\cap\Z^{n+m+1}]$ (i.e.\
$x^{[c^1,\ku{\eta}^1]} -
x^{[c^2,\ku{\eta}^2]} \in \mbox{Ker}\,\psi$), then $c^1=c^2$, and the elements
\[
\ku{\mu}^1:=\ku{\eta}^1 - \ku{\eta}(c^1),\; \ku{\mu}^2:=\ku{\eta}^2 -
\ku{\eta}(c^2) \in
C(Q)^{\veee}\cap V^\ast_{\Z}
\]
have the same image in $\C[\N^{m+1}]$.\\
In particular,
$x^{\ku{\mu}^1} - x^{\ku{\mu}^2} \in \mbox{Ker}\,(\C[C(Q)^{\veee}\cap
V^\ast_{\Z}] \rightarrow
\C[\N^{m+1}])$, and
\[
x^{[c^1,\ku{\eta}^1]} - x^{[c^2,\ku{\eta}^2]}  = x^{[c^1,\ku{\eta}(c^1)]} \cdot
\left( x^{\ku{\mu}^1} - x^{\ku{\mu}^2} \right)\, .
\vspace{2ex}
\]
\end{itemize}
\par

%%%%%%%%
% (8.4)
%%%%%%

\neu{84}
{\bf Example:}
At the end of \zitat{2}{5} we presented two decompositions of $Q_6$ into a
Minkowski sum of
lattice summands. Let us describe now the restrictions of the versal family to
the associated
components of $\bar{\CM}$:
\begin{itemize}
\item[(i)]
Putting the two triangles $R_0,R_1$ into two parallel planes contained in
$\R^3$ yields an
octahedron as the convex hull of the whole configuration. Then,
$\tilde{\sigma}$ is the
(4-dimensional) cone over this octahedron
\[
\tilde{\sigma}=\langle (0,0;1,0),\, (1,0;1,0),\,(1,1;1,0),\,(0,0;0,1),\,
(0,1;0,1),\,(1,1;0,1) \rangle\, .
\]
\item[(ii)]
Looking at the second decomposition, we have to put three line segments
$R_0,R_1,R_2$ on
three parallel 2-planes in general position inside the affine space $\R^4$.
Taking the convex
hull of this configuration yields a 4-dimensional polytope that is dual
to (triangle)$\times$(triangle).\\
Again, $\tilde{\sigma}$ is the (5-dimensional) cone over this polytope
\[
\begin{array}{r}
\tilde{\sigma}=\langle (0,0;1,0,0),\, (1,0;1,0,0),\,(0,0;0,1,0),\,(0,1;0,1,0),
\qquad\\
(0,0;0,0,1),\,(1,1;0,0,1) \rangle\, .
%\vspace{2ex}
\end{array}
\]
\end{itemize}
The total spaces over the components arise as the toric varieties defined by
$\tilde{\sigma}$. In our
example, they equal the cones over $\PP^1\times\PP^1\times\PP^1$ and
$\PP^2\times\PP^2$,
respectively.\\
\par

%%%%%%%%%%%%%%%%
%
%  Abschnitt 9
%
%%%%%%%%%%%%%%
\sect{Further examples}\label{s9}

%%%%%%%%%%
% (9.1)
%%%%%%%%

\neu{91}
Three examples of toric Gorenstein singularities arise as cones over the Del
Pezzo surfaces
obtained by blowing up $(\PP^2,{\cal O}(3))$ in one, two, or three points,
respectively. They
correspond to the following polygons:
\vspace{2ex}\\
\unitlength=0.5mm
\linethickness{0.4pt}
\begin{picture}(51.00,51.00)(-30,0)
\put(10.00,10.00){\circle*{1.00}}
\put(20.00,10.00){\circle*{1.00}}
\put(30.00,10.00){\circle*{1.00}}
\put(40.00,10.00){\circle*{1.00}}
\put(50.00,10.00){\circle*{1.00}}
\put(10.00,20.00){\circle*{1.00}}
\put(20.00,20.00){\circle*{1.00}}
\put(30.00,20.00){\circle*{1.00}}
\put(40.00,20.00){\circle*{1.00}}
\put(50.00,20.00){\circle*{1.00}}
\put(10.00,30.00){\circle*{1.00}}
\put(20.00,30.00){\circle*{1.00}}
\put(30.00,30.00){\circle*{1.00}}
\put(40.00,30.00){\circle*{1.00}}
\put(50.00,30.00){\circle*{1.00}}
\put(10.00,40.00){\circle*{1.00}}
\put(20.00,40.00){\circle*{1.00}}
\put(30.00,40.00){\circle*{1.00}}
\put(40.00,40.00){\circle*{1.00}}
\put(50.00,40.00){\circle*{1.00}}
\put(10.00,50.00){\circle*{1.00}}
\put(20.00,50.00){\circle*{1.00}}
\put(30.00,50.00){\circle*{1.00}}
\put(40.00,50.00){\circle*{1.00}}
\put(50.00,50.00){\circle*{1.00}}
\put(20.00,20.00){\line(0,1){10.00}}
\put(20.00,30.00){\line(2,1){20.00}}
\put(40.00,40.00){\line(-1,-2){10.00}}
\put(30.00,20.00){\line(-1,0){10.00}}
\put(30.00,0.00){\makebox(0,0)[cc]{Polygon $Q_4$}}
\end{picture}
\unitlength=0.5mm
\linethickness{0.4pt}
\begin{picture}(51.00,51.00)(-60,0)
\put(10.00,10.00){\circle*{1.00}}
\put(20.00,10.00){\circle*{1.00}}
\put(30.00,10.00){\circle*{1.00}}
\put(40.00,10.00){\circle*{1.00}}
\put(50.00,10.00){\circle*{1.00}}
\put(10.00,20.00){\circle*{1.00}}
\put(20.00,20.00){\circle*{1.00}}
\put(30.00,20.00){\circle*{1.00}}
\put(40.00,20.00){\circle*{1.00}}
\put(50.00,20.00){\circle*{1.00}}
\put(10.00,30.00){\circle*{1.00}}
\put(20.00,30.00){\circle*{1.00}}
\put(30.00,30.00){\circle*{1.00}}
\put(40.00,30.00){\circle*{1.00}}
\put(50.00,30.00){\circle*{1.00}}
\put(10.00,40.00){\circle*{1.00}}
\put(20.00,40.00){\circle*{1.00}}
\put(30.00,40.00){\circle*{1.00}}
\put(40.00,40.00){\circle*{1.00}}
\put(50.00,40.00){\circle*{1.00}}
\put(10.00,50.00){\circle*{1.00}}
\put(20.00,50.00){\circle*{1.00}}
\put(30.00,50.00){\circle*{1.00}}
\put(40.00,50.00){\circle*{1.00}}
\put(50.00,50.00){\circle*{1.00}}
\put(20.00,30.00){\line(0,1){10.00}}
\put(20.00,40.00){\line(1,0){10.00}}
\put(30.00,40.00){\line(1,-1){10.00}}
\put(40.00,30.00){\line(-1,-1){10.00}}
\put(30.00,20.00){\line(-1,1){10.00}}
\put(30.00,0.00){\makebox(0,0)[cc]{Polygon $Q_5$}}
\end{picture}
\unitlength=0.5mm
\linethickness{0.4pt}
\begin{picture}(51.00,51.00)(-90,0)
\put(10.00,10.00){\circle*{1.00}}
\put(20.00,10.00){\circle*{1.00}}
\put(30.00,10.00){\circle*{1.00}}
\put(40.00,10.00){\circle*{1.00}}
\put(50.00,10.00){\circle*{1.00}}
\put(10.00,20.00){\circle*{1.00}}
\put(20.00,20.00){\circle*{1.00}}
\put(30.00,20.00){\circle*{1.00}}
\put(40.00,20.00){\circle*{1.00}}
\put(50.00,20.00){\circle*{1.00}}
\put(10.00,30.00){\circle*{1.00}}
\put(20.00,30.00){\circle*{1.00}}
\put(30.00,30.00){\circle*{1.00}}
\put(40.00,30.00){\circle*{1.00}}
\put(50.00,30.00){\circle*{1.00}}
\put(10.00,40.00){\circle*{1.00}}
\put(20.00,40.00){\circle*{1.00}}
\put(30.00,40.00){\circle*{1.00}}
\put(40.00,40.00){\circle*{1.00}}
\put(50.00,40.00){\circle*{1.00}}
\put(10.00,50.00){\circle*{1.00}}
\put(20.00,50.00){\circle*{1.00}}
\put(30.00,50.00){\circle*{1.00}}
\put(40.00,50.00){\circle*{1.00}}
\put(50.00,50.00){\circle*{1.00}}
\put(20.00,20.00){\line(1,0){10.00}}
\put(30.00,20.00){\line(1,1){10.00}}
\put(40.00,30.00){\line(0,1){10.00}}
\put(40.00,40.00){\line(-1,0){10.00}}
\put(30.00,40.00){\line(-1,-1){10.00}}
\put(20.00,30.00){\line(0,-1){10.00}}
\put(30.00,0.00){\makebox(0,0)[cc]{Polygon $Q_6$}}
\end{picture}
\vspace{2ex}\\
Let us discuss these three examples:
\begin{itemize}
\item[(iv)]
The edges equal
\[
d^1=(1,0),\; d^2=(1,2),\; d^3=(-2,-1),\; d^4=(0,-1)\,,
\]
and they imply the following equations of the versal base space as closed
subscheme of
$^{\displaystyle \C^4}\!\!/\!_{\displaystyle \C\cdot (1,1,1,1)}$:
\[
%\begin{array}{cc}
t_1+t_2=2t_3,\quad t_3+t_4=2t_2,\quad t_1^2+t_2^2=2t_3^2,\quad
t_3^2+t_4^2=2t_2^2\,.
%\end{array}
\]
Using the two linear equations, only two coordinates $t:=t_1,\, \keps:=t_1-t_3$
are sufficient.
(We get the $t_i$'s back by
$t_1=t,\,t_2=t-2\keps,\,t_3=t-\keps,\,t_4=t-3\keps$). Then, the two
quadratic equations transform into $2\keps^2=0$, i.e.\ the versal base space is
a fat point.\\
On the other hand, $Q_4$ does not allow any splitting into a Minkowski sum
involving lattice
summands only. This reflects the triviality of the underlying reduced space.
(Cf.\ \zitat{9}{2}.)
\item[(v)]
The polygon $Q_5$ allows the decomposition into the sum of a triangle and a
line segment. In
particular, the reduced base space of the versal deformation of $Y_5$ has to be
a line.\\
We compute the true base space:
$d^1=(1,1),\,d^2=(-1,1),\,d^3=(-1,0),\\
d^4=(0,-1),\,d^5=(1,-1)$ yield the equations
\[
t_1-t_3=t_2-t_5=t_4-t_1\quad \mbox{and} \quad
t_1^2-t_3^2=t_2^2-t_5^2=t_4^2-t_1^2\,.
\]
With $t:=t_1,\, s_1:= t_1-t_3,\, s_2:=t_1-t_2$ and $t_1=t,\,
t_2=t-s_2,\,t_3=t-s_1,\\
t_4=t+s_1,\, t_5=t-s_1-s_2$, they turn into
\[
s_1^2=2s_1s_2=0\,.
\]
\item[(vi)]
This example was spread in the paper.
\vspace{2ex}
\end{itemize}
\par

%%%%%%%%%
% (9.2)
%%%%%%

\neu{92}
We will use the polygon $Q_4:= \mbox{Conv}\,\{(0,0),\,(1,0),\,(2,2),\,(0,1)\}$
of \zitat{9}{1}(iv) for a more detailed demonstration how the theory
works. In particular, we will describe the versal family of $Y_4$ over
$\mbox{Spec}\,^{\displaystyle \C[\keps]}\!/\!_{\displaystyle \keps^2}$:
\begin{itemize}
\item[(1)]
The $(t,\keps)$-coordinates of $V$ correspond to the linear map
\[
\left( \begin{array}{cc}
1&0\\1&-2\\1&-1\\1&-3
\end{array} \right) :
\R^2 \stackrel{\sim}{\longrightarrow} V \hookrightarrow \R^4\, .
\]
We obtain
\bean
C(Q_4) &=& \{(a,b)\in \R^2\,|\; a\geq 0,\; a-2b\geq 0,\; a-b\geq0,\; a-3b\geq
0\}\\
&=& \{(a,b)\in \R^2\,|\; a\geq 0,\; a-3b\geq 0\}\\
&=& \langle [1,0],\, [1,-3]\rangle ^{\veee} = \langle (0,-1),\,(3,1)\rangle
\subseteq \R^2\, ,
\eean
and the map $\N^4 \rightarrow C(Q_4)^{\veee}\cap V^\ast_{\Z}$ sends $e_1, e_2,
e_3, e_4$
to $[1,0]$, $[1,-2]$, $[1,-1]$, $[1,-3]$, respectively.
In particular, this map is surjective, i.e.\ $S_4=\bar{S}_4$ and
$X_4=\bar{X}_4$.
\item[(2)]
To compute the tautological cone $\tilde{C}(Q_4)$,
we need the Minkowski summands associated to the two fundamental generators of
$C(Q_4)$:
\[
(Q_4)_{(0,-1)} = \mbox{Conv}\{(0,0),\,(2,4),\,(0,3)\} ,\;
(Q_4)_{(3,1)} = \mbox{Conv}\{(0,0),\,(3,0),\,(4,2)\}\, .
\]
Hence,
\[
\begin{array}{r}
\tilde{C}(Q_4) = \langle (0,0;0,-1);\, (2,4;0,-1);\, (0,3;0,-1);\, (0,0;3,1);\,
\qquad\\
(3,0;3,1);\, (4,2;3,1) \rangle \,.
\end{array}
\]
\item[(3)]
Now, we have all information to obtain the versal family of \\
$Y_4=\mbox{Spec}\, \C[\mbox{Cone}(Q_4)^{\veee}\cap \Z^3]$:
\begin{itemize}
\item
Restrict the family $\;\mbox{Spec}\,\C[\tilde{C}(Q_4)^{\veee}\cap
\Z^4]\rightarrow
\mbox{Spec}\,\C[C(Q_4)^{\veee}\cap \Z^2]\subseteq \C^4\;$ to the subspace
$\C^2 \simeq V_{\C}\subseteq \C^4$, i.e.\ use the $(t,\keps)$-coordinates
instead of
$(t_1,t_2,t_3,t_4)$.
\item
Compose the result with the projection $\C^2 \surj \C^1\; ((t,\keps)\mapsto
\keps)$. That means,
we do no longer regard $t$ as a coordinate of the base space.
\item
Finally, we restrict our family to the closed subscheme defined by the equation
$\keps^2=0$.
\end{itemize}
\item[(4)]
To obtain equations, we could either take a closer look to the family
constructed so far, or we can
proceed more directly as described in \zitat{4}{6}, \zitat{5}{3}, and
\zitat{5}{4}:
\begin{itemize}
\item
Computing the minimal generator set of the semigroup
$\mbox{Cone}\,(Q_4)^{\veee}\cap\Z^3$,
we get the elements $[c^v;\eta_0(c^v)]$:
\[
\begin{array}{llll}
[c^1;\eta_0^1] = [0,1;0] ,&
[c^2;\eta_0^2] = [-1,1;1] ,&
[c^3;\eta_0^3] = [-2,1;2], \\
\,\![c^4;\eta_0^4] = [-1,0;2] ,&
[c^5;\eta_0^5] = [0,-1;2], &
[c^6;\eta_0^6] = [1,-2;2] ,\\
\,\![c^7;\eta_0^7] = [1,-1;1], &
[c^8;\eta_0^8] = [1,0;0] .
\end{array}
\vspace{1ex}
\]
\unitlength=0.70mm
\linethickness{0.4pt}
\begin{picture}(60.50,60.50)(-40,0)
\put(10.00,10.00){\circle*{1.00}}
\put(20.00,10.00){\circle*{1.00}}
\put(30.00,10.00){\circle*{1.00}}
\put(40.00,10.00){\circle*{1.00}}
\put(50.00,10.00){\circle*{1.00}}
\put(10.00,20.00){\circle*{1.00}}
\put(20.00,20.00){\circle*{1.00}}
\put(30.00,20.00){\circle*{1.00}}
\put(40.00,20.00){\circle*{1.00}}
\put(50.00,20.00){\circle*{1.00}}
\put(10.00,30.00){\circle*{1.00}}
\put(20.00,30.00){\circle*{1.00}}
\put(30.00,30.00){\circle*{1.00}}
\put(40.00,30.00){\circle*{1.00}}
\put(50.00,30.00){\circle*{1.00}}
\put(10.00,40.00){\circle*{1.00}}
\put(20.00,40.00){\circle*{1.00}}
\put(30.00,40.00){\circle*{1.00}}
\put(40.00,40.00){\circle*{1.00}}
\put(50.00,40.00){\circle*{1.00}}
\put(10.00,50.00){\circle*{1.00}}
\put(20.00,50.00){\circle*{1.00}}
\put(30.00,50.00){\circle*{1.00}}
\put(40.00,50.00){\circle*{1.00}}
\put(50.00,50.00){\circle*{1.00}}
\put(35.00,0.00){\makebox(0,0)[cc]{Polygon $Q_4^{\veee}$}}
\put(10.00,60.00){\circle*{1.00}}
\put(20.00,60.00){\circle*{1.00}}
\put(30.00,60.00){\circle*{1.00}}
\put(40.00,60.00){\circle*{1.00}}
\put(50.00,60.00){\circle*{1.00}}
\put(60.00,10.00){\circle*{1.00}}
\put(60.00,20.00){\circle*{1.00}}
\put(60.00,30.00){\circle*{1.00}}
\put(60.00,40.00){\circle*{1.00}}
\put(60.00,50.00){\circle*{1.00}}
\put(60.00,60.00){\circle*{1.00}}
\put(17.00,53.00){\makebox(0,0)[rb]{$z_3$}}
\put(27.00,37.00){\makebox(0,0)[rt]{$z_4$}}
\put(37.00,27.00){\makebox(0,0)[rt]{$z_5$}}
\put(53.00,17.00){\makebox(0,0)[lt]{$z_6$}}
\put(53.00,30.00){\makebox(0,0)[lc]{$z_7$}}
\put(53.00,43.00){\makebox(0,0)[lb]{$z_8$}}
\put(43.00,53.00){\makebox(0,0)[lb]{$z_1$}}
\put(30.00,53.00){\makebox(0,0)[cb]{$z_2$}}
\put(38.00,42.00){\makebox(0,0)[rb]{$t$}}
\put(20.00,50.00){\line(1,-1){30.00}}
\put(50.00,20.00){\line(0,1){20.00}}
\put(50.00,40.00){\line(-1,1){10.00}}
\put(40.00,50.00){\line(-1,0){20.00}}
\end{picture}
\vspace{2ex}
\\
Together with $[0,0;1]$, they induce coordinates $z_1,\dots,z_8,t$ on $Y_4$,
i.e.\ we have
obtained an embedding $Y_4 \hookrightarrow \C^9$.\\
(The sum of the three components of the vectors are always 1. In the figure we
have
drawn the first two coordinates.)
\item
$Y_4\subseteq \C^9$ is defined by the following 20 equations:
\[
\begin{array}{llll}
t^2-z_4z_8,&
t^2-z_1z_5,&
t^2-z_2z_7,&
z_1t-z_2z_8,\\
z_2t-z_3z_8,&
z_2t-z_1z_4,&
z_3t-z_2z_4,&
z_4t-z_3z_7,\\
z_4t-z_2z_5,&
z_5t-z_4z_7,&
z_5t-z_2z_6,&
z_6t-z_5z_7,\\
z_7t-z_5z_8,&
z_7t-z_1z_6,&
z_8t-z_1z_7,&
z_1z_3-z_2^2,\\
z_3z_5-z_4^2,&
z_4z_6-z_5^2,&
z_6z_8-z_7^2,&
z_3z_6-z_4z_5.
\end{array}
\]
\item
Choosing paths from $(0,0)\in Q_4$ to the other vertices, we obtain the list
\[
\ku{\eta}^1=[0,0,0,0], \;
\ku{\eta}^2=[1,0,0,0], \;
\ku{\eta}^3=[2,0,0,0],
\]\[
\ku{\eta}^4=[1,1,0,0]=[0,0,2,0], \;
\ku{\eta}^5=[0,2,0,0]=[0,0,1,1],
\]\[
\ku{\eta}^6=[0,0,0,2], \;
\ku{\eta}^7=[0,0,0,1], \;
\ku{\eta}^8=[0,0,0,0].
\]
\item
Now, we can lift our 20 equations to the ring
$\C[Z_1,\dots,Z_8,t_1,\dots,t_4]$:
\[
\begin{array}{llll}
t_1t_2-Z_4Z_8,&
t_2^2-Z_1Z_5,&
t_1t_4-Z_2Z_7,&
Z_1t_1-Z_2Z_8,\\
Z_2t_1-Z_3Z_8,&
Z_2t_2-Z_1Z_4,&
Z_3t_2-Z_2Z_4,&
Z_4t_3-Z_3Z_7,\\
Z_4t_2-Z_2Z_5,&
Z_5t_3-Z_4Z_7,&
Z_5t_2-Z_2Z_6,&
Z_6t_3-Z_5Z_7,\\
Z_7t_3-Z_5Z_8,&
Z_7t_4-Z_1Z_6,&
Z_8t_4-Z_1Z_7,&
Z_1Z_3-Z_2^2,\\
Z_3Z_5-Z_4^2,&
Z_4Z_6-Z_5^2,&
Z_6Z_8-Z_7^2,&
Z_3Z_6-Z_4Z_5.
\end{array}
\]
\item
Finally, we restrict the family to the versal base space
by switching to the $(t,\keps)$-coordinates
and obeying the equation $\keps^2=0$. Moreover, $t$ is no longer a coordinate
of the base space:

\[
\begin{array}{llll}
t(t-2\keps)-z_4z_8,&
t(t-4\keps)-z_1z_5,&
t(t-3\keps)-z_2z_7,\\
z_1t-z_2z_8,&
z_2t-z_3z_8,&
z_2(t-2\keps)-z_1z_4,\\
z_3(t-2\keps)-z_2z_4,&
z_4(t-\keps)-z_3z_7,&
z_4(t-2\keps)-z_2z_5,\\
z_5(t-\keps)-z_4z_7,&
z_5(t-2\keps)-z_2z_6,&
z_6(t-\keps)-z_5z_7,\\
z_7(t-\keps)-z_5z_8,&
z_7(t-3\keps)-z_1z_6,&
z_8(t-3\keps)-z_1z_7,\\
z_1z_3-z_2^2,&
z_3z_5-z_4^2,&
z_4z_6-z_5^2,\\
z_6z_8-z_7^2,&
z_3z_6-z_4z_5.
\end{array}
\]
\end{itemize}
\end{itemize}

%%%%%%%%%
% (9.3)
%%%%%%

\neu{93}
At last we want to present an example involving more than only quadratic
equations for the versal
base space. Let $Q_8$ be the ``regular'' lattice 8-gon, it is contained in two
strips of lattice
thickness 3.
\vspace{2ex}\\
\unitlength=0.50mm
\linethickness{0.4pt}
\begin{picture}(60.50,60.50)(-100,0)
\put(10.00,10.00){\circle*{1.00}}
\put(20.00,10.00){\circle*{1.00}}
\put(30.00,10.00){\circle*{1.00}}
\put(40.00,10.00){\circle*{1.00}}
\put(50.00,10.00){\circle*{1.00}}
\put(10.00,20.00){\circle*{1.00}}
\put(20.00,20.00){\circle*{1.00}}
\put(30.00,20.00){\circle*{1.00}}
\put(40.00,20.00){\circle*{1.00}}
\put(50.00,20.00){\circle*{1.00}}
\put(10.00,30.00){\circle*{1.00}}
\put(20.00,30.00){\circle*{1.00}}
\put(30.00,30.00){\circle*{1.00}}
\put(40.00,30.00){\circle*{1.00}}
\put(50.00,30.00){\circle*{1.00}}
\put(10.00,40.00){\circle*{1.00}}
\put(20.00,40.00){\circle*{1.00}}
\put(30.00,40.00){\circle*{1.00}}
\put(40.00,40.00){\circle*{1.00}}
\put(50.00,40.00){\circle*{1.00}}
\put(10.00,50.00){\circle*{1.00}}
\put(20.00,50.00){\circle*{1.00}}
\put(30.00,50.00){\circle*{1.00}}
\put(40.00,50.00){\circle*{1.00}}
\put(50.00,50.00){\circle*{1.00}}
\put(35.00,0.00){\makebox(0,0)[cc]{Polygon $Q_8$}}
\put(10.00,60.00){\circle*{1.00}}
\put(20.00,60.00){\circle*{1.00}}
\put(30.00,60.00){\circle*{1.00}}
\put(40.00,60.00){\circle*{1.00}}
\put(50.00,60.00){\circle*{1.00}}
\put(60.00,10.00){\circle*{1.00}}
\put(60.00,20.00){\circle*{1.00}}
\put(60.00,30.00){\circle*{1.00}}
\put(60.00,40.00){\circle*{1.00}}
\put(60.00,50.00){\circle*{1.00}}
\put(60.00,60.00){\circle*{1.00}}
\put(20.00,40.00){\line(1,1){10.00}}
\put(30.00,50.00){\line(1,0){10.00}}
\put(40.00,50.00){\line(1,-1){10.00}}
\put(50.00,40.00){\line(0,-1){10.00}}
\put(50.00,30.00){\line(-1,-1){10.00}}
\put(40.00,20.00){\line(-1,0){10.00}}
\put(30.00,20.00){\line(-1,1){10.00}}
\put(20.00,30.00){\line(0,1){10.00}}
\end{picture}
\vspace{2ex}\\
$Q_8$ admits three maximal Minkowski decompositions into a sum of lattice
summands:
\vspace{-6ex}
\begin{itemize}
\item[(i)]
$Q_8\;=$
\unitlength=0.50mm
\linethickness{0.4pt}
\begin{picture}(40.50,40.50)(6,23)
\put(10.00,10.00){\circle*{1.00}}
\put(20.00,10.00){\circle*{1.00}}
\put(30.00,10.00){\circle*{1.00}}
\put(40.00,10.00){\circle*{1.00}}
\put(10.00,20.00){\circle*{1.00}}
\put(20.00,20.00){\circle*{1.00}}
\put(30.00,20.00){\circle*{1.00}}
\put(40.00,20.00){\circle*{1.00}}
\put(10.00,30.00){\circle*{1.00}}
\put(20.00,30.00){\circle*{1.00}}
\put(30.00,30.00){\circle*{1.00}}
\put(40.00,30.00){\circle*{1.00}}
\put(10.00,40.00){\circle*{1.00}}
\put(20.00,40.00){\circle*{1.00}}
\put(30.00,40.00){\circle*{1.00}}
\put(40.00,40.00){\circle*{1.00}}
\put(20.00,20.00){\line(0,1){10.00}}
\end{picture}
+
\unitlength=0.50mm
\linethickness{0.4pt}
\begin{picture}(40.50,40.50)(6,23)
\put(10.00,10.00){\circle*{1.00}}
\put(20.00,10.00){\circle*{1.00}}
\put(30.00,10.00){\circle*{1.00}}
\put(40.00,10.00){\circle*{1.00}}
\put(10.00,20.00){\circle*{1.00}}
\put(20.00,20.00){\circle*{1.00}}
\put(30.00,20.00){\circle*{1.00}}
\put(40.00,20.00){\circle*{1.00}}
\put(10.00,30.00){\circle*{1.00}}
\put(20.00,30.00){\circle*{1.00}}
\put(30.00,30.00){\circle*{1.00}}
\put(40.00,30.00){\circle*{1.00}}
\put(10.00,40.00){\circle*{1.00}}
\put(20.00,40.00){\circle*{1.00}}
\put(30.00,40.00){\circle*{1.00}}
\put(40.00,40.00){\circle*{1.00}}
\put(30.00,20.00){\line(-1,0){10.00}}
\end{picture}
+
\unitlength=0.50mm
\linethickness{0.4pt}
\begin{picture}(40.50,40.50)(6,23)
\put(10.00,10.00){\circle*{1.00}}
\put(20.00,10.00){\circle*{1.00}}
\put(30.00,10.00){\circle*{1.00}}
\put(40.00,10.00){\circle*{1.00}}
\put(10.00,20.00){\circle*{1.00}}
\put(20.00,20.00){\circle*{1.00}}
\put(30.00,20.00){\circle*{1.00}}
\put(40.00,20.00){\circle*{1.00}}
\put(10.00,30.00){\circle*{1.00}}
\put(20.00,30.00){\circle*{1.00}}
\put(30.00,30.00){\circle*{1.00}}
\put(40.00,30.00){\circle*{1.00}}
\put(10.00,40.00){\circle*{1.00}}
\put(20.00,40.00){\circle*{1.00}}
\put(30.00,40.00){\circle*{1.00}}
\put(40.00,40.00){\circle*{1.00}}
\put(30.00,20.00){\line(-1,1){10.00}}
\end{picture}
+
\unitlength=0.50mm
\linethickness{0.4pt}
\begin{picture}(40.50,40.50)(6,23)
\put(10.00,10.00){\circle*{1.00}}
\put(20.00,10.00){\circle*{1.00}}
\put(30.00,10.00){\circle*{1.00}}
\put(40.00,10.00){\circle*{1.00}}
\put(10.00,20.00){\circle*{1.00}}
\put(20.00,20.00){\circle*{1.00}}
\put(30.00,20.00){\circle*{1.00}}
\put(40.00,20.00){\circle*{1.00}}
\put(10.00,30.00){\circle*{1.00}}
\put(20.00,30.00){\circle*{1.00}}
\put(30.00,30.00){\circle*{1.00}}
\put(40.00,30.00){\circle*{1.00}}
\put(10.00,40.00){\circle*{1.00}}
\put(20.00,40.00){\circle*{1.00}}
\put(30.00,40.00){\circle*{1.00}}
\put(40.00,40.00){\circle*{1.00}}
\put(20.00,20.00){\line(1,1){10.00}}
\end{picture}
\item[(ii)]
$Q_8\;=$
\unitlength=0.50mm
\linethickness{0.4pt}
\begin{picture}(40.50,40.50)(6,23)
\put(10.00,10.00){\circle*{1.00}}
\put(20.00,10.00){\circle*{1.00}}
\put(30.00,10.00){\circle*{1.00}}
\put(40.00,10.00){\circle*{1.00}}
\put(10.00,20.00){\circle*{1.00}}
\put(20.00,20.00){\circle*{1.00}}
\put(30.00,20.00){\circle*{1.00}}
\put(40.00,20.00){\circle*{1.00}}
\put(10.00,30.00){\circle*{1.00}}
\put(20.00,30.00){\circle*{1.00}}
\put(30.00,30.00){\circle*{1.00}}
\put(40.00,30.00){\circle*{1.00}}
\put(10.00,40.00){\circle*{1.00}}
\put(20.00,40.00){\circle*{1.00}}
\put(30.00,40.00){\circle*{1.00}}
\put(40.00,40.00){\circle*{1.00}}
\put(20.00,20.00){\line(1,1){10.00}}
\put(30.00,30.00){\line(0,-1){10.00}}
\put(30.00,20.00){\line(-1,0){10.00}}
\end{picture}
+
\unitlength=0.50mm
\linethickness{0.4pt}
\begin{picture}(40.50,40.50)(6,23)
\put(10.00,10.00){\circle*{1.00}}
\put(20.00,10.00){\circle*{1.00}}
\put(30.00,10.00){\circle*{1.00}}
\put(40.00,10.00){\circle*{1.00}}
\put(10.00,20.00){\circle*{1.00}}
\put(20.00,20.00){\circle*{1.00}}
\put(30.00,20.00){\circle*{1.00}}
\put(40.00,20.00){\circle*{1.00}}
\put(10.00,30.00){\circle*{1.00}}
\put(20.00,30.00){\circle*{1.00}}
\put(30.00,30.00){\circle*{1.00}}
\put(40.00,30.00){\circle*{1.00}}
\put(10.00,40.00){\circle*{1.00}}
\put(20.00,40.00){\circle*{1.00}}
\put(30.00,40.00){\circle*{1.00}}
\put(40.00,40.00){\circle*{1.00}}
\put(20.00,20.00){\line(1,1){10.00}}
\put(30.00,30.00){\line(-1,0){10.00}}
\put(20.00,30.00){\line(0,-1){10.00}}
\end{picture}
+
\unitlength=0.50mm
\linethickness{0.4pt}
\begin{picture}(40.50,40.50)(6,23)
\put(10.00,10.00){\circle*{1.00}}
\put(20.00,10.00){\circle*{1.00}}
\put(30.00,10.00){\circle*{1.00}}
\put(40.00,10.00){\circle*{1.00}}
\put(10.00,20.00){\circle*{1.00}}
\put(20.00,20.00){\circle*{1.00}}
\put(30.00,20.00){\circle*{1.00}}
\put(40.00,20.00){\circle*{1.00}}
\put(10.00,30.00){\circle*{1.00}}
\put(20.00,30.00){\circle*{1.00}}
\put(30.00,30.00){\circle*{1.00}}
\put(40.00,30.00){\circle*{1.00}}
\put(10.00,40.00){\circle*{1.00}}
\put(20.00,40.00){\circle*{1.00}}
\put(30.00,40.00){\circle*{1.00}}
\put(40.00,40.00){\circle*{1.00}}
\put(30.00,20.00){\line(-1,1){10.00}}
\end{picture}
\item[(iii)]
$Q_8\;=$
\unitlength=0.50mm
\linethickness{0.4pt}
\begin{picture}(40.50,40.50)(6,23)
\put(10.00,10.00){\circle*{1.00}}
\put(20.00,10.00){\circle*{1.00}}
\put(30.00,10.00){\circle*{1.00}}
\put(40.00,10.00){\circle*{1.00}}
\put(10.00,20.00){\circle*{1.00}}
\put(20.00,20.00){\circle*{1.00}}
\put(30.00,20.00){\circle*{1.00}}
\put(40.00,20.00){\circle*{1.00}}
\put(10.00,30.00){\circle*{1.00}}
\put(20.00,30.00){\circle*{1.00}}
\put(30.00,30.00){\circle*{1.00}}
\put(40.00,30.00){\circle*{1.00}}
\put(10.00,40.00){\circle*{1.00}}
\put(20.00,40.00){\circle*{1.00}}
\put(30.00,40.00){\circle*{1.00}}
\put(40.00,40.00){\circle*{1.00}}
\put(20.00,30.00){\line(1,0){10.00}}
\put(30.00,30.00){\line(0,-1){10.00}}
\put(30.00,20.00){\line(-1,1){10.00}}
\end{picture}
+
\unitlength=0.50mm
\linethickness{0.4pt}
\begin{picture}(40.50,40.50)(6,23)
\put(10.00,10.00){\circle*{1.00}}
\put(20.00,10.00){\circle*{1.00}}
\put(30.00,10.00){\circle*{1.00}}
\put(40.00,10.00){\circle*{1.00}}
\put(10.00,20.00){\circle*{1.00}}
\put(20.00,20.00){\circle*{1.00}}
\put(30.00,20.00){\circle*{1.00}}
\put(40.00,20.00){\circle*{1.00}}
\put(10.00,30.00){\circle*{1.00}}
\put(20.00,30.00){\circle*{1.00}}
\put(30.00,30.00){\circle*{1.00}}
\put(40.00,30.00){\circle*{1.00}}
\put(10.00,40.00){\circle*{1.00}}
\put(20.00,40.00){\circle*{1.00}}
\put(30.00,40.00){\circle*{1.00}}
\put(40.00,40.00){\circle*{1.00}}
\put(30.00,20.00){\line(-1,0){10.00}}
\put(20.00,20.00){\line(0,1){10.00}}
\put(20.00,30.00){\line(1,-1){10.00}}
\end{picture}
+
\unitlength=0.50mm
\linethickness{0.4pt}
\begin{picture}(40.50,40.50)(6,23)
\put(10.00,10.00){\circle*{1.00}}
\put(20.00,10.00){\circle*{1.00}}
\put(30.00,10.00){\circle*{1.00}}
\put(40.00,10.00){\circle*{1.00}}
\put(10.00,20.00){\circle*{1.00}}
\put(20.00,20.00){\circle*{1.00}}
\put(30.00,20.00){\circle*{1.00}}
\put(40.00,20.00){\circle*{1.00}}
\put(10.00,30.00){\circle*{1.00}}
\put(20.00,30.00){\circle*{1.00}}
\put(30.00,30.00){\circle*{1.00}}
\put(40.00,30.00){\circle*{1.00}}
\put(10.00,40.00){\circle*{1.00}}
\put(20.00,40.00){\circle*{1.00}}
\put(30.00,40.00){\circle*{1.00}}
\put(40.00,40.00){\circle*{1.00}}
\put(20.00,20.00){\line(1,1){10.00}}
\end{picture}
\vspace{6ex}
\pagebreak
\end{itemize}
The decompositions (i), (ii) and (i), (iii) are refinements of the coarser
decompositions
\vspace{-5ex}
\\
$Q_8\;=$
\unitlength=0.50mm
\linethickness{0.4pt}
\begin{picture}(50.50,45.50)(6,23)
\put(10.00,5.00){\circle*{1.00}}
\put(20.00,5.00){\circle*{1.00}}
\put(30.00,5.00){\circle*{1.00}}
\put(40.00,5.00){\circle*{1.00}}
\put(50.00,5.00){\circle*{1.00}}
\put(10.00,15.00){\circle*{1.00}}
\put(20.00,15.00){\circle*{1.00}}
\put(30.00,15.00){\circle*{1.00}}
\put(40.00,15.00){\circle*{1.00}}
\put(50.00,15.00){\circle*{1.00}}
\put(10.00,25.00){\circle*{1.00}}
\put(20.00,25.00){\circle*{1.00}}
\put(30.00,25.00){\circle*{1.00}}
\put(40.00,25.00){\circle*{1.00}}
\put(50.00,25.00){\circle*{1.00}}
\put(10.00,35.00){\circle*{1.00}}
\put(20.00,35.00){\circle*{1.00}}
\put(30.00,35.00){\circle*{1.00}}
\put(40.00,35.00){\circle*{1.00}}
\put(50.00,35.00){\circle*{1.00}}
\put(10.00,45.00){\circle*{1.00}}
\put(20.00,45.00){\circle*{1.00}}
\put(30.00,45.00){\circle*{1.00}}
\put(40.00,45.00){\circle*{1.00}}
\put(50.00,45.00){\circle*{1.00}}
\put(20.00,15.00){\line(0,1){10.00}}
\put(20.00,25.00){\line(1,1){10.00}}
\put(30.00,35.00){\line(1,0){10.00}}
\put(40.00,35.00){\line(0,-1){10.00}}
\put(40.00,25.00){\line(-1,-1){10.00}}
\put(30.00,15.00){\line(-1,0){10.00}}
\end{picture}
+
\unitlength=0.50mm
\linethickness{0.4pt}
\begin{picture}(40.50,40.50)(6,23)
\put(10.00,10.00){\circle*{1.00}}
\put(20.00,10.00){\circle*{1.00}}
\put(30.00,10.00){\circle*{1.00}}
\put(40.00,10.00){\circle*{1.00}}
\put(10.00,20.00){\circle*{1.00}}
\put(20.00,20.00){\circle*{1.00}}
\put(30.00,20.00){\circle*{1.00}}
\put(40.00,20.00){\circle*{1.00}}
\put(10.00,30.00){\circle*{1.00}}
\put(20.00,30.00){\circle*{1.00}}
\put(30.00,30.00){\circle*{1.00}}
\put(40.00,30.00){\circle*{1.00}}
\put(10.00,40.00){\circle*{1.00}}
\put(20.00,40.00){\circle*{1.00}}
\put(30.00,40.00){\circle*{1.00}}
\put(40.00,40.00){\circle*{1.00}}
\put(30.00,20.00){\line(-1,1){10.00}}
\end{picture}
and
$\quad Q_8\;=$
\unitlength=0.50mm
\linethickness{0.4pt}
\begin{picture}(50.50,45.50)(6,23)
\put(10.00,5.00){\circle*{1.00}}
\put(20.00,5.00){\circle*{1.00}}
\put(30.00,5.00){\circle*{1.00}}
\put(40.00,5.00){\circle*{1.00}}
\put(50.00,5.00){\circle*{1.00}}
\put(10.00,15.00){\circle*{1.00}}
\put(20.00,15.00){\circle*{1.00}}
\put(30.00,15.00){\circle*{1.00}}
\put(40.00,15.00){\circle*{1.00}}
\put(50.00,15.00){\circle*{1.00}}
\put(10.00,25.00){\circle*{1.00}}
\put(20.00,25.00){\circle*{1.00}}
\put(30.00,25.00){\circle*{1.00}}
\put(40.00,25.00){\circle*{1.00}}
\put(50.00,25.00){\circle*{1.00}}
\put(10.00,35.00){\circle*{1.00}}
\put(20.00,35.00){\circle*{1.00}}
\put(30.00,35.00){\circle*{1.00}}
\put(40.00,35.00){\circle*{1.00}}
\put(50.00,35.00){\circle*{1.00}}
\put(10.00,45.00){\circle*{1.00}}
\put(20.00,45.00){\circle*{1.00}}
\put(30.00,45.00){\circle*{1.00}}
\put(40.00,45.00){\circle*{1.00}}
\put(50.00,45.00){\circle*{1.00}}
\put(20.00,35.00){\line(1,0){10.00}}
\put(30.00,35.00){\line(1,-1){10.00}}
\put(40.00,25.00){\line(0,-1){10.00}}
\put(40.00,15.00){\line(-1,0){10.00}}
\put(30.00,15.00){\line(-1,1){10.00}}
\put(20.00,25.00){\line(0,1){10.00}}
\end{picture}
+
\unitlength=0.50mm
\linethickness{0.4pt}
\begin{picture}(40.50,40.50)(6,23)
\put(10.00,10.00){\circle*{1.00}}
\put(20.00,10.00){\circle*{1.00}}
\put(30.00,10.00){\circle*{1.00}}
\put(40.00,10.00){\circle*{1.00}}
\put(10.00,20.00){\circle*{1.00}}
\put(20.00,20.00){\circle*{1.00}}
\put(30.00,20.00){\circle*{1.00}}
\put(40.00,20.00){\circle*{1.00}}
\put(10.00,30.00){\circle*{1.00}}
\put(20.00,30.00){\circle*{1.00}}
\put(30.00,30.00){\circle*{1.00}}
\put(40.00,30.00){\circle*{1.00}}
\put(10.00,40.00){\circle*{1.00}}
\put(20.00,40.00){\circle*{1.00}}
\put(30.00,40.00){\circle*{1.00}}
\put(40.00,40.00){\circle*{1.00}}
\put(20.00,20.00){\line(1,1){10.00}}
\end{picture},
\vspace{5ex}
\\
respectively.\\
\par
These facts translate directly into the geometry of the reduced base space of
the versal deformation
of $Q_8$:
\begin{itemize}
\item
It is embedded in some affine space $\C^5$ and equals the union of a 3-plane
with two
2-planes (through $0\in \C^5$).
\item
The two 2-planes each have a common line with the 3-dimensional component.
However, they
intersect each other in $0\in \C^5$ only.
\vspace{2ex}
\end{itemize}
On the other hand, we can write down the equations of the true versal base
space (as a closed
subscheme of $^{\displaystyle \C^8}\!\!/\!_{\displaystyle \C\cdot
(1,\dots,1)}$):
\[
t_1^k +t_2^k+t_8^k=t_4^k+t_5^k+t_6^k,\quad
t_2^k+t_3^k+t_4^k=t_6^k+t_7^k+t_8^k\quad
(k=1,2,3)\,.
\]
\\
\par

%%%%%%
%
% References
%
%%%%%%%

%}
\end{document}